\documentclass[a4paper,12pt]{amsart}\usepackage{amsfonts,amssymb,bbm}
\usepackage{enumerate,mathrsfs,xfrac,xspace,hyperref} \UseRawInputEncoding

\def\C{\Bbb{C}}
\def\k{\Bbbk}

\def\bk{{\bar{\k}}}\def\N{\mathbb{N}}
\def\bxi{{\bar{\xi}}}\def\txi{{\tilde{\xi}}}

\def\R{\Bbb{R}}

\def\di{\partial}

\def\bl{\langle}\def\br{\rangle}

\def\RmX{R^{\oplus m}_X}

\newcommand{\quots}[2]{{\footnotesize\left.\raisebox{0.4ex}{$#1$}\! / \!\raisebox{-0.4ex}{$#2$}\right.}}

\def\tA{\tilde{A}}\def\tB{\tilde{B}}\def\tf{{\tilde{f}}}\def\tF{\tilde{F}}
\def\tJ{{\tilde{J}}}
 \def\tn{{\tilde{n}}}
\def\bPhi{{\bar{\Phi}}}

\def\tv{\tilde{v}}\def\tx{{\tilde{x}}}\def\tX{{\tilde{X}}}\def\ty{{\tilde{y}}}\def\tY{{\tilde{Y}}}
\def\tz{{\tilde{z}}}

\def\hg{{\hat{g}}}
\def\hk{{\hat\k}}\def\hQ{\hat{Q}}\def\hR{{\widehat{R}}}
\def\hU{\hat{U}}\def\hV{\hat{V}}
\def\hx{\hat{x}}\def\hX{{\hat{X}}}\def\htX{{\hat{\tX}}}\def\hy{\hat{y}}\def\hY{{\hat{Y}}}\def\htY{{\hat{\tY}}}\def\hz{\hat{z}}
\def\hPhi{\hat{\Phi}}

\def\be{\beta}\def\de{\delta}\def\De{\Delta}\def\Om{\Omega}
\def\ep{\epsilon}

\def\ca{\mathfrak a}
\def\cC{\mathscr C}

\def\cK{{\mathscr K}\!}\def\cL{\mathscr L}\def\cR{\mathscr{R}}

\def\cm{{\frak m}}
\def\cU{\mathcal U}

\def\ud{{\underline{d}}}

\def\um{{\underline{m}}}
\def\uv{\underline{v}}

\newcommand{\ber}{\begin{array}{l}}\newcommand{\eer}{\end{array}}
\newcommand{\bpm}{\begin{pmatrix}}\newcommand{\epm}{\end{pmatrix}}
\newcommand{\bbm}{\begin{bmatrix}}\newcommand{\ebm}{\end{bmatrix}}
\newcommand{\bM}{\begin{matrix}}\newcommand{\eM}{\end{matrix}}
\newcommand{\bee}{\begin{enumerate}}\newcommand{\eee}{\end{enumerate}}
\newcommand{\bei}{\begin{itemize}}\newcommand{\eei}{\end{itemize}}

\def\wrt{with respect to }
\def\sset{\subset}\def\sseteq{\subseteq}\def\ssetneq{\subsetneq}\def\smin{\setminus}

\def\Maps{Maps (X,Y)}
\def\MapX{Maps\big(X,(\k^m,o)\big)}
\def\Mapk{Maps\big((\k^n,o),(\k^m,o)\big)}
\def\RAP{$\cR$.AP\xspace}\def\RAPP{$\cR$.APP\xspace}\def\LAP{$\cL$.AP\xspace}\def\LAPP{$\cL$.APP\xspace}
\def\LRAP{$\cL\cR$.AP\xspace}\def\LRAPP{$\cL\cR$.APP\xspace}

\newtheorem{Lemma}{Lemma}[section]\newcommand{\bel}{\begin{Lemma}}\newcommand{\eel}{\end{Lemma}}
\newtheorem{Theorem}[Lemma]{Theorem}\newcommand{\bthe}{\begin{Theorem}}\newcommand{\ethe}{\end{Theorem}}
\newtheorem{Proposition}[Lemma]{Proposition}\newcommand{\bprop}{\begin{Proposition}}\newcommand{\eprop}{\end{Proposition}}
\newtheorem{Corollary}[Lemma]{Corollary}\newcommand{\bcor}{\begin{Corollary}}\newcommand{\ecor}{\end{Corollary}}
\newtheorem{Definition}[Lemma]{Definition}\newcommand{\bed}{\begin{Definition}}\newcommand{\eed}{\end{Definition}}
\newtheorem{Definition-Proposition}[Lemma]{Definition-Proposition}

\def\bpr{{\em\noindent Proof.\ }}
\newcommand{\epr}{{\hfill\ensuremath\blacksquare}}
\newtheorem{Remark}[Lemma]{Remark}\newcommand{\beR}{\begin{Remark}\rm}\newcommand{\eeR}{\end{Remark}}
\newtheorem{Example}[Lemma]{Example}\newcommand{\bex}{\begin{Example}\rm}\newcommand{\eex}{\end{Example}}
\newtheorem{Problem}[Lemma]{Problem}\newcommand{\bprob}{\begin{Problem}\rm}\newcommand{\eprob}{\end{Problem}}

\newcommand{\bet}{\begin{tabular}{lll}}\newcommand{\eet}{\end{tabular}}
\newcommand{\beq}{\vspace{-0.1cm}\begin{equation}}\newcommand{\eeq}{\vspace{-0.1cm}\end{equation}}

 \renewcommand{\stackrel}[2]{\ \lower 0.4ex \hbox{$\mathrel{\mathop{#2}\limits^{\scriptscriptstyle {#1}}}$}\ }
\newcommand{\isom}[1]{\hbox{$\xrightarrow[\,\smash{\raisebox{1.15ex}{\ensuremath{\scriptstyle\sim}}}\,]{#1}$}}

\title[]{R\MakeLowercase{esults on left-right} A\MakeLowercase{rtin approximation  for    algebraic morphisms\\and for  analytic morphisms of weakly-finite singularity type}}
\author[]{D\MakeLowercase{mitry} K\MakeLowercase{erner}}
\address{Department of Mathematics, Ben Gurion University of the Negev, P.O.B. 653, Be'er Sheva 84105, Israel. dmitry.kerner@gmail.com}
\date{\today\ \  filename: \jobname.tex}
\thanks{I was supported by the Israel Science Foundation, grants No.  1910/18 and 1405/22}


\begin{document}
\begin{abstract}
 The classical  Artin approximation (AP) reads: any formal solution of a system of (analytic, resp. algebraic) equations of implicit function type is approximated by ``ordinary"  solutions (i.e. analytic, resp. algebraic).

  Morphisms of scheme-germs, e.g.  $\Mapk,$ are usually studied up to the left-right equivalence. The natural question  is the left-right version of Artin approximation: when is the formal left-right equivalence of morphisms   approximated by the ``ordinary"  (i.e. analytic, resp. algebraic) equivalence?
 In this case the standard Artin approximation is not directly  applicable, as the involved (functional) equations are
  not of implicit function type.   Moreover, the naive extension does not hold in the  analytic case, because of Osgood-Gabrielov-Shiota examples.

  The left-right version of Artin approximation  (\LRAP)  was established by M. Shiota for  morphisms that are either Nash or [real-analytic and of finite singularity type].

\medskip

    We establish \LRAP\ and its stronger version of P\l oski (\LRAPP) for   $\Maps,$ where $X,Y$ are  analytic/algebraic germs of schemes of any characteristic.
    More precisely:
   \bei
   \item
   \LRAP, \LRAPP, the inverse Artin approximation (and its P\l oski's version) hold for    algebraic morphisms and for finite analytic morphisms.
     \item
      \LRAP holds for analytic morphisms of weakly-finite singularity type. (For $char>0$ we impose certain integrability condition.)
\eei 

\medskip

This latter class of morphisms of   ``weakly-finite singularity type" (which we introduce) is of separate importance.
 It extends naturally the traditional class of morphisms of ``finite singularity type", while preserving their non-pathological behaviour.
 The definition goes via the higher critical loci and higher discriminants of morphisms with singular targets.

 We establish basic properties of these critical loci.
 In particular: any map   is  finitely (right) determined by its higher critical loci.
\end{abstract}
\maketitle \vspace{-0.8cm}
\setcounter{secnumdepth}{6} \setcounter{tocdepth}{1}\tableofcontents
\vspace{-0.7cm}

\section{Introduction}\label{Sec.Introduction}
\subsection{The ordinary AP and APP}\label{Sec.Introduction.ordinary.AP.APP} Consider a system of implicit function equations, $F(x,y)=0.$ Here $x,y$ are multi-variables, while $F$ is a (finite) vector of power series.
 Assume $F$ is either analytic or algebraic, i.e. its entries belong to one of the rings:
\bee[\!\!$\bullet$\!]
\item $\k\{x,y\},$ the  convergent power series ($\k$ is a normed field, complete \wrt its norm)
\item $\k\bl x,y\br,$ the algebraic power series ($\k$ is any field). For $\k=\R$ these are the Nash power series.
\eee
The classical approximation theorems \cite{Artin.68},\cite{Artin.69} read:  any formal solution, $F(x,\hy(x))\!=\!0,$ is
 approximated ($x$-adically) by ordinary solutions, $y(x)$ (analytic or algebraic power series).
 (See \S\ref{Sec.Background.AP.classical} for more detail.)

 Moreover, there exists a ``parameterized" (analytic/algebraic) solution, $F(x,y(x,z))=0,$ whose specialization is the initial formal solution. Namely, there exists a formal power series $\hz(x)\in \k[\![x]\!]$ satisfying:  $y(x,\hz(x))=\hy(x).$
 This says (roughly): the formal solutions are not ``intrinsic"\!\!, they are just ``unfortunate" specialization of the ordinary (analytic, resp. algebraic) solutions.
 See \cite{Ploski1974} for $\k\{x,y\}$ and \cite{Popescu.85}, \cite{Popescu.86} for $\k\bl x,y\br.$
  (See also \cite{Ogoma.94}, \cite{Swan.98}.)

\medskip

 The subject ``Artin approximations" has been studied extensively, see e.g. \cite{Popescu}, \cite{Rond} for numerous results and further references. It became an everyday tool in Local Algebraic/Analytic Geometry and in Commutative Algebra.

\vspace{-0.2cm}

\subsection{The inverse AP}\label{Sec.Intro.Inverse.AP} Fix some (analytic, resp. algebraic) power series, $f(x)\!=\!f_1(x),\!\dots\!,f_m(x).$
 One way to check their dependence is to resolve the functional  equation $\Phi(f(x))\!=\!0,$ with the unknown $\Phi.$  The ``inverse Artin question", \!\cite[pg.13-08]{Grothendieck.61}, \!asks: is any formal solution, $\hPhi(f(x))\!=\!0,$ \!approximated by ordinary solutions, \!$\Phi(f(x))\!=\!0?$
 (Here $\Phi$ is analytic, resp. algebraic.)

 The condition $ \Phi(f(x))=0$ is not an  implicit function equation on $\Phi$.
 Thus the inverse AP  is not a direct consequence of the classical version of AP.

 This inverse Artin approximation holds in the algebraic case, the proof is based on the nested AP of \cite{Popescu.85}, \cite{Popescu.86}.
  In the analytic case the inverse AP holds for the ring $\C\{x\}$ of dimension$\le\!\!2$, \cite{Abhyankar-van der Put.70}. It
  fails in $dim\ge3$, see \cite{Gabrielov.71}, \cite{Gabrielov.73}.

\medskip

The inverse AP was studied in several works,  see the references in \cite{Izumi.07},
 and is still under investigation, see   \cite{Bel.Cur.Ron.}, \cite{Alonso.Castro.Hauser.Koutschan}.
  For linear equations over admissible rings (in particular $\k\bl x\br,$ $\k\{x\}$) this inverse AP problem is transformed into a nested AP problem,  \cite{Cas-Jim.Pop.Ron}.

  It seems, no P\l oski version of the inverse AP is known.

 \vspace{-0.1cm}
\subsection{The left-right AP and APP}\label{Sec.Introduction.Left.Right.AP}
The germs of maps, $\Mapk,$ are traditionally studied up to the left-right equivalence,
 $f\sim \Phi_Y\circ f\circ\Phi_X^{-1},$ \cite{AGLV}, \cite{Mond-Nuno}.  Here $f$ is a vector of (analytic, resp. algebraic) power series, while $\Phi_X\in Aut_{(\k^n,o)},\Phi_Y\in Aut_{(\k^m,o)}$
  are (analytic, resp. algebraic) automorphisms, i.e. coordinate changes in the source and in the target.
  The basic task is to verify whether two maps are equivalent, $f\stackrel{\cL\cR}{\sim}\tf,$ or to bring a map to a particular ``normal" form.
  The arguments are often inductive, they provide only the formal equivalence.
  Hence the fundamental question is the  Artin approximation   and its Artin-P\l oski-Popescu version, as in \S\ref{Sec.Introduction.ordinary.AP.APP},  for the equation $\Phi_Y\circ \tf=f\circ\Phi_X$:

\bei
\item {\pmb\LRAP:}  {\em Is every formal solution, $\hPhi_Y\circ \tf= f\circ\hPhi_X$,    approximated by  ordinary solutions?}
\item {\pmb\LRAPP:} {\em Is \!every formal solution a specialization of a parameterized ordinary solution?}
\eei
Here the maps $f,\tf$ are prescribed, while the morphisms $\Phi_X,\Phi_Y$ are unknowns.
 We do not require   $\Phi_X,\Phi_Y$  to be isomorphisms. But if $\hPhi_X,\hPhi_Y$ are invertible, then so are their approximations.

\medskip

As in \S\ref{Sec.Intro.Inverse.AP},
 the condition $ \Phi_Y \circ  \tf =    f \circ  \Phi_X$ is not an  implicit function equation.
 Thus the properties \LRAP, \LRAPP are not direct consequences of the classical version of AP.

\beR
\bee[\bf i.]
\item
Occasionally we add the condition $\Phi_X = Id_X,$ resp.  $\Phi_Y=Id_Y,$  to get the property \LAP\!\!, resp. \RAP\!\!\!.
 Geometrically the property \RAP addresses isomorphisms of scheme-germs over a given base: any formal isomorphism
 $[X$ over $Y]\isom{}[\tX$ over $Y]$ is approximated by ordinary isomorphisms.
 For \LRAP one allows the base change: any formal isomorphism  $[X$ over $Y]\isom{}[\tX$ over $\tY]$ is approximated by ordinary isomorphisms.

\item   The classical Artin approximation is related to the \RAP question: when is the formal right equivalence of maps, $f \circ \hPhi^{-1}_X = \tf,$ approximated by  ordinary equivalence? (See Example \ref{Ex.RAP}.)

    Similarly,  the  inverse Artin approximation can be stated as the   \LAP  question:  when  is
 the formal left equivalence of maps,  $\hPhi_Y \circ  f = \tf,$ approximated  by   ordinary  equivalence? (See  Remark  \ref{Rem.LRAP}.)

\item Various approximation questions become trivial for finitely determined maps. A map $f$ is called finitely-$\cL\cR$-determined if the equality modulo higher order terms, $f\equiv \tf\ mod(x)^N$ (for $N\gg1$), implies the equivalence, $f\stackrel{\cL\cR}{\sim}\tf.$
    In other words, the group orbit $\cL\cR(f)$ contains an $(x)$-adically open neighborhood of $f$ inside $\Mapk.$

For these maps the property \LRAP holds trivially. See \cite[Chapter 6]{Mond-Nuno} for the classical determinacy (with $\k=\R,\C$), see
 \cite[\S7]{Kerner.Group.Orbits} for the case of arbitrary $\k.$
\eee\eeR

\medskip

 The property \LRAP was studied by M. Shiota for $Maps\big((\R^n,o),(\R^m,o)\big).$
\bei 
\item For  Nash-germs, $\R\bl x\br,$ the property \LRAP holds for any map, \cite[Theorem 4]{Shiota.2010}.

Though stated for $\R\bl x\br$, the proof is actually characteristic-free and works over any field.
\item
In the $\R$-analytic case, $\R\{x\},$ the property \LRAP does not always hold, \cite[Fact.1.4]{Shiota.1998}, \cite[pg. 1061]{Shiota.2010}.
 The counterexample is based on the counterexample of \cite{Gabrielov.71} to the Inverse Artin question, \S\ref{Sec.Intro.Inverse.AP}.

 However, the property \LRAP holds for maps of finite singularity type\footnote{\label{note1}i.e. $f$ is finite \wrt   contact equivalence,
       i.e. the complexified subgerm $V(f)\sset (\C^n,o)$ is either an ICIS or   zero-dimensional, see \S\ref{Sec.Weak.Fin.Finite.Sing.Type}.
        An $\cL\cR$-finitely-determined map is necessarily of finite singularity type. But the converse does not hold in ``most" cases, see e.g. Chapter 7 of \cite{Mond-Nuno}.}, \cite[Fact.1.7]{Shiota.1998}. That proof, \cite[pg.119-126]{Shiota.1998}, is heavily $\R$-analytic.
 \eei 

 \medskip

 The global version of \LRAP, for $\Maps$ with $X,Y$ compact Nash/real-analytic manifolds, is  more complicated and is
  also treated in \cite{Shiota.1998}, \cite{Shiota.2010}.

  These approximation results were extended  to equations of type
   $F(x,\!\Phi(x),\!\Psi\!\circ\! g(x))\!=\!0$ in \cite{Fichou.Quarez.Shiota}. Here $F,g$ are prescribed $\R$-Nash power series, while  $\Phi,\Psi$ are unknowns.

 \vspace{-0.1cm}
\subsection{Results}\label{Sec.Introduction.Results}To study left-right equivalence of morphisms $\Maps$, \cite{Kerner.Group.Orbits}, and their deformations, \cite{Kerner.Unfoldings},   we  need the relevant Artin-P\l oski-Popescu  approximations over any \!$\k$.

  \vspace{-0.1cm}
  \subsubsection{\!\!Extending the results of \cite{Shiota.1998},\cite{Shiota.2010} in several directions}  We study $\Maps$
 instead of $\Mapk.$ Here $X\!:=\!Spec(R_X)$ and $Y\!:=Spec(R_Y)$ are scheme-germs with arbitrary singularities.
  Algebraically the maps are local homomorphisms of $\k$-algebras, $Hom_\k(R_Y,R_X).$

\bee[\!\!$\bullet$\!]
\item   The case of algebraic power series, $R_X\!\!=\!\!\quots{\k\bl x\br}{J_X}\!,$ $R_Y\!\!=\!\!\quots{\k\bl y\br}{J_Y}\!,$
     with $\k$   a field or an excellent Henselian local ring over a field.
\\{\bf Theorem \ref{Thm.LRAP.for.Nash}.}  {\em
    The properties \LRAP, \LRAPP, \LAP,   \LAPP, the inverse Artin approximation (and its P\l oski version) hold  for   algebraic morphisms $ \Maps.$}

\item   The finite analytic case, $R_X\!\!=\!\!\quots{\k\{x\}}{J_X}\!,$ $R_Y\!\!=\!\!\quots{\k\{y\}}{J_Y}\!,$\! with  $\k$    a  complete  normed field or    a Henselian DVR over a field.
\\{\bf Theorem \ref{Thm.LRAP.finite.analytic.maps}.}
{\em  The properties \LRAP, \LRAPP, \LAP, \LAPP, the inverse Artin approximation (and its P\l oski version)  hold for   finite  analytic morphisms,   $ \Maps.$}

 \item  (The non-finite analytic case)
\\{\bf Theorem \ref{Thm.LRAP.weakly.finite.analytic}.} {\em  The property \LRAP (with $\hPhi_Y,\hPhi_X$ isomorphisms)   holds   for analytic morphisms $ \Maps$ of weakly-finite singularity type.}
\\
(For $char>0$ we impose additional conditions: ``approximate integrability" of certain derivations, and ``approximate lifting" of certain automorphisms.)
\eee

\medskip

\noindent The natural questions now are about \!\LRAPP, the inverse Artin approximation, \!and \!its P\l oski \!version, for morphisms \!of weakly-finite  singularity \!type. These demand additional ingredients.
 \!(Left for \!the \!subsequent work.)

  \vspace{-0.1cm}
  \subsubsection{The weakly-finite singularity type}\label{Sec.Intro.Weakly.Finite.Sing.Type} It is the natural extension  of the classical ``finite singularity type", seems not previously studied.
  We sketch now the ideas of the definition.  The main ingredient is the critical locus for morphisms with {\em non-smooth targets}.
   It is defined via the truncated cotangent complex.
\bee[\!\!$\bullet$\!]
\item  The critical module of a dominant morphism $f:X\to Y$ is defined via the exact sequence $ Der_X \stackrel{f}{\to} T_{X\to Y}\to \cC_{X,f}\to0.$
     Here $Der_X:=Der_\k(R_X)$ is the module of derivations, while $T_{X\to Y}$ is the tangent space to $\Maps$ at $f.$
\\
 The critical locus, $Crit_{X}f\sseteq X,$  is the support of the critical module, $\cC_{X,f}.$ For a $\k$-smooth target, $Y=(\k^m,o),$ one gets the classical critical locus, defined by the Jacobian ideal of the map $f.$ Otherwise $Crit_X$ is usually smaller than the branch locus (the points where the module of relative differentials, $ \Om_{X/Y},$ is non-free) and than  the tangential branch locus
  (the points where the module of relative vector fields, $T_{X/Y},$ is non-free).

 The pathology $Crit_X f=X$ can occur for non-finite morphisms due to Osgood-type examples, i.e. dominant maps $X\to Y$ with $dim(Y)>dim(X).$
  For $char>0$ the pathology $Crit_X f=X$ can occur even for finite morphisms.

 The discriminant of $f$ is the image of the critical locus, $\De:=f(Crit_{X}f)\sset Y.$

\item  Take the restricted map   $f_|:Crit_{X}f\to \De.$ One gets the ``second" critical locus and its discriminant,
 $Crit_{X}f\supseteq Crit_{Crit,f_|}\to \De_{f_|}\sseteq \De.$ Iterating this we get  the higher critical modules,  $\cC_j$, the higher critical loci,
  $X=Crit_0\supseteq Crit_1\supseteq \cdots,$  and the higher discriminants,
  $Y:=\De_0\supseteq \De_1\supseteq\cdots .$

 \item A map $f:X\to Y$  is called {\em of weakly-finite singularity type} (w.f.s.t.) if its restriction $f_|:Crit_r\to \De_r$ is finite for some $r\ge0.$
\bee[\bf-]
\item   For $r=0$ one gets just a finite map.
\item For $r=1$ one gets a map of finite singularity type (f.s.t.), extending the classical f.s.t. notion, \cite[pg.224]{Mond-Nuno}.
\item For scalar-valued maps, $Maps(X,(\k^1,o)),$ the   notions f.s.t. and w.f.s.t. coincide.
\item For higher dimensional targets  the condition w.f.s.t.  is   much weaker  than f.s.t.
\eee
 Thus, even in the classical case, $\Mapk,$ for $\k=\R,\C$, Theorem \ref{Thm.LRAP.weakly.finite.analytic} is stronger than the known results.

\eee

The notion of $Crit_{X}f$ for morphisms with singular target is known, \cite{Kallstrom}. But we could not find a reference for even the basic properties.
 We establish these in \S\ref{Sec.Weak.Fin.Critical.Locus}-\S\ref{Sec.Weak.Fin.Weakly.Finite.Sing.Type}.

  \vspace{-0.2cm}
\subsubsection{}\label{Sec.Introduction.Auxil.Results} Two auxiliary results in this paper are of separate importance.
\bee[\!\!$\bullet$\!]
\item (Approximate lifting of automorphisms.) Take a subgerm $Z=Spec(\quots{R_X}{I})\sset X.$ If $X$ is singular, then not all automorphisms of $Z$  extend to $X.$ (There is no lifting from the group of automorphisms $Aut_Z$ to $Aut_X.$)
    Even for $N$-th thickening, $Z_N:=Spec(\quots{R_X}{I^N})\sset X,$ with $N\gg1,$ the lifting from $Aut_{Z_N}$ to $Aut_X$ does not need to exist.
 But an approximation of each automorphism can be lifted:
\\
{\bf Lemma \ref{Thm.Lifting.Property.Automorphisms}.} (for $char\!=\!0$) {\em
For any ideal $I\!\sset \!R$ and a suitable $N\!\gg\!1,$ any automorphism $\bPhi\!\in\! Aut_\k(\quots{R}{I^N})$
 has an approximate lifting. Namely, there exists $\Phi\!\in\! Aut_\k(R)$ satisfying:  $\Phi(I)\!=\!I$ and \!$[\Phi]\!=\! [\bPhi]\!\in\! Aut_\k(\quots{R}{I}).$}
\\
This property would follow (in any characteristic) by the Strong Artin Approximation for the filtration $I^\bullet\!\sset\! R.$ But that version of SAP does not hold,
 \cite[Example 3.20.]{Rond}.

 \medskip

\item  Any morphism is finitely-right-determined by its higher critical loci.
\\
{\bf Lemma \ref{Thm.Higher.Critical.Locus.Fin.Determ.Map}.} (roughly) {\em
 Suppose    dominant maps $X\stackrel{f,\tf}{\to}Y$ have coinciding higher critical loci and discriminants, i.e. (for a certain $j\ge0$):
 \[Crit_i(f)\!=\!Crit_i(\tf)\quad  \text{ and }\quad \De_i(f)\!=\!\De_i(\tf) \quad \text{ for }\quad i=0,\dots,j.
 \]
 If $f\!-\!\tf\!\in\! Fitt_0(\cC_{j+1})^N\!\cdot\! \RmX$ for some $N\gg1,$ then $\tf$ is $\cR$-equivalent to $f,$ i.e. $\tf\!=\!f\!\circ\!\Phi_X$ for $\Phi_X\!\in\! Aut_X.$}\hfill  (Here $\cC_{j+1}$ is the higher critical module, as in \S\ref{Sec.Intro.Weakly.Finite.Sing.Type}.)
 \\
For  $j\!=\!0$ this  is a well-known fact: ``any morphism is finitely-$\cR$-determined by its critical locus", \cite[Example 4.2]{Kerner.Group.Orbits}.
\eee

\subsection{A bit about the (purely algebraic) proofs of the approximation results}
\bee[\!\!$\bullet$\!]
\item (The proof of Theorem \ref{Thm.LRAP.for.Nash}) We reformulate the $\cL\cR$-condition, $\Phi_Y\circ f=\tf\circ\Phi_X,$
 as a   homomorphism of algebras with base change. This translates the initial  functional equation  into a problem of implicit function type,   with nested solutions. Then, for the algebraic case ($\quots{\k\bl x\br}{J_X},\quots{\k\bl y\br}{J_Y}$), one uses the nested Artin approximation/nested P\l oski theorem.

\item (The proof of Theorem \ref{Thm.LRAP.finite.analytic.maps})
In the   analytic case, ($\quots{\k\{ x\}}{J_X},\quots{\k\{ y\}}{J_Y}$), the nested Artin approximation does not hold, \cite[\S5.3]{Rond}.
 However, we can consider the map $R_X\stackrel{\Phi_X}{\to}R_\tX$ not as a homomorphism of algebras, but as  a homomorphism of (finitely-generated) modules,
  $mod$-$R_Y\ni R_X\to R_\tX\in mod$-$R_\tY,$
  with a base change $R_Y\to R_\tY.$ This translates the initial conditions into a system of implicit function equations.
 Eventually  this homomorphism of modules  is promoted to a homomorphism of algebras, ensuring $f \sim \tf .$

\item (The proof of Theorem \ref{Thm.LRAP.weakly.finite.analytic}) By the classical Artin approximation we can assume: all the higher critical loci of $f,\tf$ coincide. And similarly for the higher discriminants. Thus we get two finite analytic maps $f_|,\tf_|: Crit_r\to \De_r,$ that are formally left-right equivalent.
    By  Theorem \ref{Thm.LRAP.finite.analytic.maps} these restrictions are analytically left-right equivalent,
      via the group $Aut_{Crit_r}\times Aut_{\De_r}$.
    But this equivalence is not directly liftable to $Aut_X\times Aut_Y,$ see \S\ref{Sec.Introduction.Auxil.Results}.

    Therefore we start from $N$-th thickening of $Crit_r$, and repeat the argument. Then Lemma \ref{Thm.Lifting.Property.Automorphisms} ensures an approximate lifting
 to $Aut_X\times Aut_Y.$  Eventually one gets:     $f-\tf\in (I_{Crit_r})^N,$ for $N\gg1.$ And now invoke Lemma  \ref{Thm.Higher.Critical.Locus.Fin.Determ.Map}.
 \eee

\subsection{Applications}
 Singularities of mapping-germs were classically studied in the real/complex cases. One of the main tools was vector field integration.
 The case of arbitrary base field (of arbitrary characteristic) has been under investigation in the last 20 years, see e.g. Chapter 3 of \cite{Gr.Lo.Sh}. In the absence of vector field integration, the proofs were by induction,  thus only for rings of formal power series.

  Therefore the left-right Artin approximation (in all its versions) is of key importance. E.g. in \cite{Kerner.Unfoldings} we develop the unfolding theory of (analytic/algebraic) map-germs, extending the classical criteria of triviality/versality/stability  from $\R,\C$ to an arbitrary field $\k$. The proofs there are heavily based on Theorems \ref{Thm.LRAP.for.Nash}, \ref{Thm.LRAP.finite.analytic.maps}, \ref{Thm.LRAP.weakly.finite.analytic}.

\subsection{Acknowledgements} Thanks are due to A.-F. Boix for valuable knowledge and remarks, to G. Rond for the nice survey \cite{Rond}, and to D. Popescu for important advices.

We thank the referee for the very careful reading  and numerous helpful remarks.

 Special thanks are to the organizers of
 the ``\L ojasiewicz Workshop"  (December 2022, Gda\'{n}ńsk-Krak\'{o}ów-{\L}Ł\'{o}ód\'{z}ź-Warszawa),
  during which the essential advance was done, and to the Alibaba Kebab Shop (Marsza\l kowska 83).

\section{Definitions and preliminaries}
\subsection{Notation  and conventions (Germs of schemes and their morphisms)}\label{Sec.Background.Schemes.Maps}
 For the general introduction to  $\Mapk$
 in the real/complex-analytic case see \cite{Mond-Nuno}.  See \cite{Kerner.Group.Orbits} for some definitions and results over arbitrary fields.

 \bee[\bf i.]
\item We use   multivariables, $x=(x_1,\dots,x_n),$  $y=(y_1,\dots,y_m),$  $F=(f_1,\dots,f_r).$
 The notation  $F(x,y)\in (x,y)^N,$ resp. $y(x)\in (x)^N,$ means: every component of the vector $F$ (resp. $y(x)$) satisfies this condition.

\noindent  Denote the maximal ideal of $\k$  by $\cm_\k.$ (If $\k$ is a field, then $\cm_\k=0$.) We always assume $F(x,y)\in \cm_\k+(x,y).$

\item Consider (analytic, resp. algebraic) morphisms of scheme-germs, $\Maps$. Here:
\bei
\item  $X\!=\!Spec(R_X)$, where  $R_X\!=\!\quots{\k\{x\}}{J_X},$ resp.  $R_X\!=\!\quots{\k\bl x\br}{J_X}.$
 (We assume $J_X\!\sseteq\! (x)^2$.)
\item  $Y=Spec(R_Y)$, where $R_Y\!=\!\quots{\k\{y\}}{J_Y},$ resp.   $R_Y\!=\!\quots{\k\bl y\br}{J_Y}.$
 (We assume $J_Y\!\sseteq\! (y)^2$.)

\item
For analytic power series, $\quots{\k\{x\}}{J_X},$  $\k$ is either a normed field  (complete \wrt its norm), or a Henselian DVR over a field.
 (The norm is always non-trivial.)

 The \!main relevant examples of \!DVR are  the rings
 $\k_o[\![\tau]\!], \k_o\bl \tau\br, \k_o\{\tau\}.$ (Here  $\k_o$ is a field.)

   \item For algebraic power series, $\quots{\k\bl x\br}{J_X},$   $\k$ is any  field or an excellent Henselian local ring over a field.
 See \cite{Moret-Bailly} for more detail. Occasionally we call this case ``$\k$-Nash" to avoid any confusion with the algebraic rings like  $\k[x]$ or  $\k[x]_{(x)}.$
\eei

\medskip

  In the $\k$-smooth case,  $J_X=0,J_Y=0,$ we denote  $(\k^n,o):=Spec(R_X)$, $(\k^m,o):=Spec(R_Y).$ These are formal/analytic/$\k$-Nash germs,
     depending on the context.

  Otherwise the singularities of $X,Y$ are arbitrary.

Recall that the rings $R:=\k[\![x]\!],\k\{x\},\k\bl x\br$ admit substitutions: if $f\in R$ and $g_1,\dots,g_n\in (x)\sset R,$ then $f(g_1,\dots,g_n)\in R.$

  \item
A (formal/analytic/$\k$-Nash) map $f\in \Maps$ is defined by the corresponding homomorphism of $\k$-algebras, $f^\sharp\in Hom_\k(R_Y,R_X).$ All our homomorphisms are local.
\\In the   case   $J_X=0,J_Y=0 $    one gets the maps of $\k$-smooth germs,  $\Mapk.$

 Fix   embeddings $X\sseteq (\k^n,o),$ $Y\sseteq (\k^m,o)$ and some coordinates $(x_1,\dots,x_n)$ on $(\k^n,o),$ and $(y_1,\dots,y_m)$ on  $(\k^m,o).$
 They are sent to the generators of ideals in $R_X,$ $R_Y$, abusing notations we write just  $(x)=(x_1,\dots,x_n)\sset R_X$ $(y)=(y_1,\dots,y_m)\sset R_Y.$
     If $\k$ is a field, then $(x),(y)$ are the maximal ideals. Otherwise the maximal ideals are $\cm_X:=\cm_\k+(x)\sset R_X$
      and $\cm_Y:=\cm_\k+(y)\sset R_Y.$  When only $R_X$ is involved we write just $\cm\sset R_X.$

  Then a map $f:X\to Y$ is represented by
  the tuple of power series, $(f_1,\dots,f_m)\in \cm\cdot\RmX.$
 This embeds the space of maps as a subset $\Maps\sseteq \cm\cdot\RmX.$
 Explicitly, for $R_Y=\quots{\k[\![y]\!]}{J_Y},\quots{\k\{y\}}{J_Y},\quots{\k\bl y\br}{J_Y},$ one has:
 \beq
 \Maps=Hom_\k(R_Y,R_X)=\{f\in \cm\cdot \RmX|\ f^\sharp(J_Y)\sseteq J_X.\}
 \eeq
 If the target is $\k$-smooth, $Y=(\k^m,o),$ (i.e. $R_Y=\k[\![y]\!],$ resp. $\k\{y\},$ resp. $\k\bl y\br$),
  then we get the $R_X$-module   $\MapX= \cm\cdot\RmX.$

\item
 The  (formal/analytic/$\k$-Nash) automorphisms of $X$ over $Spec(\k)$ are defined via their algebraic counterparts, $\k$-linear automorphisms of the local ring, $Aut_X:=Aut_\k(R_X).$
 Similarly we denote  $Aut_Y:=Aut_\k(R_Y).$
  In the $\k$-smooth case, $\Mapk,$ these automorphisms are the coordinate changes in the source and the target.

These automorphisms act on the space of maps, defining the right and left group actions:
\bei
\item $\cR:=Aut_{X}: =Aut_\k(R_X)\circlearrowright \Maps$, by $(\Phi_X,f)\to f\circ \Phi^{-1}_X$;
\item $\cL:=Aut_{Y}: =Aut_\k(R_Y)\circlearrowright \Maps$, by $ (\Phi_Y,f)\to \Phi_Y\circ f$.
\eei
These define the right, left, and left-right equivalence of maps, $\cL\times\cR\!: f\sim \Phi_Y\circ f\circ \Phi^{-1}_X.$

      \eee

\subsection{The classical versions of Artin approximation}\label{Sec.Background.AP.classical} Below we assume $\hy(x)\in \cm\cdot\hk[\![x]\!].$
\bee[\bf i.]
\item ({\em The ordinary version for analytic power series,} \cite[Theorem 1.1]{Denef-Lipshitz})
 Let $\k$ be a complete normed field or a Henselian DVR.
 Consider equations $F(x,y)=0$ for $F(x,y)\in (x,y)\sset\k\{x,y\}.$
 Then any formal solution, $\hy(x)\in \cm\cdot \hk[\![x]\!],$ is $\cm$-adically approximated by  analytic solutions.

\item ({\em The nested version for algebraic power series,} \cite{K.P.P.R.M.}, \cite{Popescu.86}, \cite[Prop.3.5]{Popescu}, \cite[Theorem 5.8]{Rond})
Let $\k$ be a field or
 an excellent Henselian local ring.
 Consider equations $F(x,y)=0$ for $F(x,y)\in (x,y)\sset\k\bl x,y\br.$  Suppose these equations have a  formal  solution  $\hy(x)\in \cm$
 that is nested:
  \beq\label{Eq.nested.solution.AP}
 \hy_1(x)\in \hk[\![x_1,\dots,x_{i_1}]\!], \dots, \hy_m(x)\in \hk[\![x_1,\dots,x_{i_m}]\!], \text{ for some } 1\le i_1\le \cdots\le i_m\le n.
 \eeq
   This formal solution    is $\cm$-adically  approximated by  nested algebraic solutions,
    $y_1(x)\in \k\bl x_1,\dots,x_{i_1}\br,$ \dots, $y_r(x)\in \k\bl x_1,\dots,x_{i_m}\br.$

\item ({\em The strong Artin approximation,} \cite{Pfister-Popescu}, \cite[Theorem 7.1]{Denef-Lipshitz})

Let $\k$ be a field or a complete DVR, and take the power series
 $F(x,y)\in \k[\![x,y]\!].$ There exists a function $\be_F:\N\to \N$ satisfying: for every $d\in \N$ any approximate solution, $F(x,\bar{y}(x))\in \cm^{\be_F(d)},$ is $\cm^\bullet$-approximated by an ordinary solution, i.e. $F(x,y(x))=0$ with $y(x)-\bar{y}(x)\in \cm^d.$

 \item ({\em Parametrization of solutions}, \cite{Ploski1974}, \cite{Popescu.85}, \cite{Popescu.86}) Given a formal solution, $F(x,\hy(x))=0,$ $\hy(x)\in \cm\cdot \k[\![x]\!],$  there exist power series $y(x,z)\in \k\{x,z\}^{\oplus m},$ resp. $\k\bl x,z\br^{\oplus m},$
      and (formal) power series $\hz(x)\in\cm\cdot \k[\![x]\!]^{\oplus c},$ satisfying:
     \beq\label{Eq.Ploski.theorem}
     F(x,y(x,z))=0,\quad \text{ and }\quad  \hy(x)=y(x,\hz(x)) .
     \eeq

     \noindent({\em The nested version for $\k\bl x \br$},   Theorem 11.4 of \cite{Spivakovsky.99}, Theorem 2.1 of \cite{Bilski.Parusiński.Rond})
     If the formal solution $\hy(x)$ is nested, as in \eqref{Eq.nested.solution.AP}, then the parameterizations $y(x,z),$ $ \hz(x)$
     in equation \eqref{Eq.Ploski.theorem} are nested as well.
     Namely:
     \bei
     \item     $
      y_1(x,z)\in \k\bl x_1,\dots,x_{i_1},z_1,\dots,z_{j_1}\br,$ \dots, $
      y_m(x,z)\in \k\bl x_1,\dots,x_{i_m},z_1,\dots,z_{j_m}\br,$ \quad and
      \item
      $\hz_1(x)\dots \hz_{j_1}(x)\in \k[\![x_1,\dots,x_{i_1}]\!],$
       $\hz_{j_1+1}(x)\dots \hz_{j_2}(x)\in \k[\![x_1,\dots,x_{i_2}]\!],$ \dots,
       \eei
       for some $1\le i_1\le \cdots\le i_m\le m$ and $1\le j_1\le \cdots\le j_m\le c.$
\eee

\beR\label{Rem.Background.Embedding.as.IF.eq}
\bee[\!\!\!\!\bf i.\!]

\item Let $R_X$ be one of the rings $\quots{\k[\![x]\!]}{J_X},\quots{\k\{x\}}{J_X},\quots{\k\bl x\br}{J_X}.$ Consider the systems of implicit function equations, $F(x,y)=0.$ Here $F\in R_X[\![y]\!],$ resp. $R_X\{ y\}, R_X\bl y\br.$
    When $J_X\neq0$ this system of equations is understood in the following sense.
     Take any representative $\tF$ of $F$ over $\k[\![x,y]\!],$ resp. $\k\{x,y\} ,\k\bl x,y\br.$
      Then $F(x,y)=0$ means: $\tF(x,y)\equiv0$ mod $J_X,$ i.e. $\tF(x,y)\in J_X.$

\item
 (Presenting the condition $\Phi(I)\sseteq J$ as an implicit function equation)

 Denote by $R$ the rings $\quots{\k[\![x]\!]}{\ca},$   $\quots{\k\bl x\br}{\ca},$ $\quots{\k\{x\}}{\ca}.$
  Fix two ideals $I,J\sset R$
  and  a homomorphism $\Phi:R\to R.$ For the fixed generators, $x=(x_1,\dots,x_n)$ in $R,$ we present $\Phi$ by power series,
   $\{\Phi(x_i)=\tx_i(x)\}_i.$  Fix some (finite set of) generators, $\{q^I_i\}$ of $I,$ and $\{q^J_j\}$ of $J.$
    Then the condition $\Phi(I)\sseteq J$ transforms into the set of equations:
   \beq
  [ \quad \Phi(q^I_i(x))=\joinrel= ]\hspace{1cm}q^I_i(\tx(x))=\sum \tz_{ji}\cdot q^J_j(x),\quad \forall\ i.
   \eeq
This is a system of implicit function equations for the unknowns $\{\tx_i\},\{\tz_{ji}\}$ in $R.$
 \eee

 Hence the Artin approximation is applicable also to the condition $\Phi(I)\sseteq J.$
\eeR

\bex\label{Ex.RAP} ({\em The property ``\RAP\!\!"}) Let $X,\tX$ be analytic (resp. Nash) scheme-germs, \S\ref{Sec.Background.Schemes.Maps}.ii. Any formal morphism $\hat\Phi: \hat\tX\to \hX$ (i.e. a homomorphism of the $\cm$-adic completions $\hat\Phi^\sharp: \hR_X\to \hR_\tX$) is approximated by  ordinary morphisms, $\tX\to X$ (i.e.  homomorphisms $ \Phi^\sharp:  R_X\to R_\tX,$ $x\to x(\tx)$).
 Namely, for each $d$ there exists an analytic (resp. Nash) morphism $\Phi_d:\tX\to X$ satisfying $\hat\Phi^\sharp(x)-\Phi^\sharp_d(x)\in (x)^d.$

 In particular (for $X=\tX$), any formal automorphism, $\hat\Phi\in Aut_\hk(\hR_X),$ is approximated by ordinary automorphisms.

 ({\em The property ``\RAPP\!\!"}) Moreover, there exists a parameterized family of morphisms, $\Psi: \tX\times(\k^c_z,o)\to X,$
  i.e. $\Psi^\sharp: R_X\to R_\tX\{z\}$ (resp. $\Psi^\sharp: R_X\to R_\tX\bl z\br$) satisfying: $\hPhi^\sharp(x) =\Psi(x,\hz(x)),$
   for a specialization $\hz(x)\in \hR_X.$
\eex

\subsection{The condition $f\circ\Phi_X=\Phi_Y\circ \tf$ as a morphism of algebras with base-change}\label{Sec.Background.Commut.Diagram}
The map $f\in \Maps$ induces on $R_X$ the structure of  $R_Y$-algebra.
 The structure map is $R_Y\ni g(y)\stackrel{f^\sharp}{\to} g(f)\in R_X$.
 We observe: $f\stackrel{\cL\cR}{\sim}\tf$ iff the corresponding algebras admit an isomorphism with (invertible) base-change, as on the diagram
\beq\label{Eq.diagram}
\bM f\in &R_X&\isom{x\to x(\tx)}& R_\tX&\ni \tf
\\\uparrow\quad &\uparrow f^\sharp&  & \tf^\sharp\uparrow&\quad \uparrow
\\ y\in &R_Y&\isom{y\to y(\ty)}& R_\tY&\ni \ty\eM
 \quad .
\eeq
Indeed, the commutativity of the diagram means $f(x(\tx))=y(\tf(\tx))$, i.e.:
   $f(\Phi_X(\tx))=\Phi_Y(\tf(\tx)),$ where $\Phi_X(\tx)=x(\tx)$ and $\Phi_Y(\ty)=y(\ty).$

\

In some cases one imposes the condition $f\circ\Phi_X=\Phi_Y\circ \tf$ without invertibility of $\Phi_X,\Phi_Y.$ Then this diagram expresses just a  homomorphism of algebras with base-change.

\section{Left-right Artin approximation in the Nash case and in the finite analytic case}\label{Sec.LRAP}

\bed\label{Def.LRAP} The left-right Artin approximation, \LRAP, holds  for a map $f\in \Maps$ if for any map $\tf\in \Maps,$
 and the equation $ \Phi_Y\circ \tf=f\circ  \Phi_X,$ any formal solution $\hat\Phi_Y\circ \tf=f\circ \hat\Phi_X$
   is $\cm$-adically approximated by   ordinary solutions.
 \eed
\noindent\mbox{Namely, for any $d\in \N$ there exist  endomorphisms of $\k$-algebras $\Phi_X\circlearrowright R_X$ and $\Phi_Y\circlearrowright R_Y$,
  satisfying:}
  \beq
  \Phi_X(x)-\hat\Phi_X(x)\in \cm_X^d, \quad\quad\quad\quad\quad \Phi_Y(y)-\hat\Phi_Y(y)\in \cm_Y^d, \quad\quad\quad\quad\quad
    \Phi_Y\circ \tf=f\circ  \Phi_X.
\eeq
\bed\label{Def.LRAPP}
The left-right parameterized approximation, \LRAPP,
  holds  for a map $f\in \Maps,$ if for any map $\tf\in \Maps $
 and the equation $ \Phi_Y\circ \tf=f\circ  \Phi_X,$ any formal solution $\hat\Phi_Y\!\circ \!\tf\!=\!f\!\circ\! \hat\Phi_X$
  is a specialization of a parameterized ordinary solution, $\Phi_{X,z_X}$, $\Phi_{Y,z_Y}.$
   \eed
Namely, there exist Nash homomorphisms  $\Phi^\sharp_{Y,z_Y}:R_Y\to R_\tY\bl z_Y\br ,$
 $\Phi^\sharp_{X,z_X}:R_X\to R_\tX\bl z_X\br,$    (resp. their analytic versions), acting by $y\to y(\ty,z_Y) $ and $x\to x(\tx,z_X),$ and
  formal power series $\hz_X(x)\in \hR_X,$ $\hz_Y(y)\in \hR_Y,$ that altogether satisfy:
\beq
\Phi_{Y,z_Y}\circ \tf=f\circ  \Phi_{X,z_X},  \quad\quad\quad
\hPhi_X=\Phi_{X,\hz_X(x)},\quad\quad\quad
\hPhi_Y=\Phi_{Y,\hz_Y(y)}.
\eeq
(Here $f,\tf,\hPhi_X,\hPhi_Y$ are prescribed, while $\Phi_{Y,z_Y},\Phi_{X,z_X},\hz_X(x),\hz_Y(y)$ are unknowns.)

\beR\label{Rem.LRAP}
\bee[\!\!\!\!\bf i.\!]
\item
We do not impose invertibility of $\Phi_X,\Phi_Y$.
But if $\hat\Phi_X,\hat\Phi_Y$ are invertible, then their approximations $\Phi_X,\Phi_Y$ are ring isomorphisms/automorphisms,
    $\Phi_X\!\in\! \cR$ and $\Phi_Y\!\in\! \cL.$

  \item Imposing $\Phi_Y=Id_Y,$ one gets
  the properties \RAP, \RAPP. They are implied by AP, APP, Example \ref{Ex.RAP}. Conversely, AP can be stated as a condition of nested \RAP-type: \\suppose a map $F(x,y)$ is formally $\cR$-equivalent (by $y\to y-\hg(x)$) to a map of type $\sum y_i\cdot \tF_i(x,y).$ Is this (formal, nested) equivalence approximated by ordinary nested equivalence?

\item Imposing $\Phi_X=Id_X,$ one gets   the properties \LAP, \LAPP.
\\
Let $Y=(\k^m,o)$. The property \LAP   implies the inverse Artin approximation for $f.$ The property \LAPP implies the  P\l oski version of
the inverse  Artin approximation.
\\\bpr Take a formal relation $\hQ(f)=0.$ If $\hQ(y)\sseteq (y)^2$ then define $\hat\Phi_Y:=Id+\hQ\in \hat\cL.$ Thus $\hat\Phi_Y(f)=f.$ Approximate $\hat\Phi$ by an ordinary automorphism, $\Phi_Y(f)=f.$ Then $Q:=\Phi_Y-Id_Y$ is the needed approximation of $\hQ.$

     If $\hQ(y)\not\sseteq (y)^2$ then apply the previous argument to   $ \{y_i \cdot \hQ(y)\} .$

The P\l oski version: using \LAPP one gets $\Phi_{Y,z}$ satisfying: $\Phi_{Y,z}(f)=f$ and $\Phi_{Y,\hz(y)}=\hPhi_Y.$ Then take $Q_{Y,z}:=\Phi_{Y,z}-Id_Y.$
     \epr
\item
Suppose for a given class of maps the condition $\Phi_Y\circ \tf=f\circ\Phi_X $ can be translated into a (finite) system of implicit function equations, without nested conditions. Then \LRAP, \LRAPP hold for this class of  maps.

Moreover, the properties \RAP and \LAP  hold as well.
 Indeed, to the condition $\Phi_Y \circ  \tf\! = \!f  \circ   \Phi_X$ one just appends the
   condition $\Phi_Y(y) \!= \!y$ (resp. $\Phi_X(x)\!\!=\!\!x$). And the later conditions are implicit function equations.
\item
Definition \ref{Def.LRAP}  addresses approximation \wrt the filtrations $\cm^\bullet_X\sset R_X, $ $\cm^\bullet_Y\sset R_Y.$
 One can obtain better approximations, for stronger filtrations, as follows. Suppose $f\circ\hPhi_X=\hPhi_Y\circ\tf,$ where the formal automorphisms satisfy:
 \[
 \hPhi_X\equiv Id_X\ mod\ I_X\ \rm{and}\  \hPhi_Y\equiv Id_Y\ mod\ I_Y, \quad \rm{for} \quad \sqrt{I_X}\ssetneq\cm_X \ \rm{and} \  \sqrt{I_Y}\ssetneq\cm_Y.
 \]
If \LRAP   holds for $f,$ for the filtrations $\cm^\bullet_X,$ $\cm^\bullet_Y,$ then it holds for the filtrations
   $I_X\cdot \cm^\bullet_X,$ $I_Y\cdot \cm^\bullet_Y.$ The proof follows that of \cite{Belitskii.Boix.Kerner}.
   The detail will be given in our next work.
\eee
\eeR

\subsection{Approximations for algebraic power series, $R_X\!=\!\quots{\k\bl x\br}{J_X}$ and $R_Y=\quots{\k\bl y\br}{J_Y}$}\label{Sec.LRAP.for.k<x>}

\bthe\label{Thm.LRAP.for.Nash}
 The properties \LRAP, \LRAPP, \LAP,   \LAPP, the inverse Artin approximation (and its P\l oski version) hold  for any Nash map $f\in \Maps.$
 \ethe
\bpr
In view of \S\ref{Sec.Background.Commut.Diagram} we present $R_X$ as an $R_Y$-algebra, $R_X=\quots{R_Y\bl x\br}{(J_X+(y-f(x)))}.$
 Similarly $R_\tX=\quots{R_\tY\bl \tx\br}{(J_\tX+(\ty-\tf(\tx)))}.$
 We have formal horizontal homomorphisms  in the diagram \eqref{Eq.diagram},
 and  want to approximate these by   ordinary homomorphisms of the  algebras $R_X,R_\tX$ with a base change.
 The needed homomorphism is a map $\psi:R_X\to R_\tX$ presentable by (algebraic) power series, $\psi(y)=y(\ty)$ and $\psi(x)=x(\tx).$
  And $\psi$ must satisfy:
  \beq\label{Eq.inside.proof.LRAP.Nash.Case}
  \psi(J_Y)\sseteq J_\tY,\hspace{1cm} \psi(J_X+(y-f(x)))\sseteq J_\tX+(\ty-\tf(\tx)).
  \eeq
 These inclusions of ideals are implicit function equations, see    Remark \ref{Rem.Background.Embedding.as.IF.eq}.ii.
 Indeed:
 \bei
 \item[-] for $J_X=0$ and $J_Y=0$    the condition is: $y(\ty)-f(x(\tx))=s(\tx,\ty)\cdot (\ty-\tf(\tx)),$ with the unknowns
  $y(\ty), x(\tx),s(\tx,\ty).$
\item[-] for  $J_X\neq0$ or $J_Y\neq0$   take the additional conditions and further auxiliary variables.
\eei
 By our assumption these implicit function equations  admit a formal solution,
\beq
\hat\psi(y)=\hy(\ty)\in \k[\![\ty]\!],\quad\quad \hat\psi(x)=\hx(\tx)\in \k[\![\tx]\!], \quad\quad \hz(\tx,\ty)\in \k[\![\tx,\ty]\!].
\eeq
    We can think of it as a formal nested solution,  $\hat\psi(y)=\hy(\ty)$ and $\hat\psi(x)=\hx(\tx,\ty)\in \k[\![\tx,\ty]\!].$

Altogether, we have transformed the initial condition into a system of implicit function equations with nested solution.
\bee[\!\!$\bullet$\!]
\item (\LRAP)
 Apply the nested Artin approximation, \S\ref{Sec.Background.AP.classical}.ii.,
   to approximate this formal nested solution by an algebraic nested solution, $y(\ty)\!\in\! \k\bl \ty\br$ and
  $x(\tx,\ty)\!\in\! \k\bl \tx,\ty\br.$ Explicitly:
  \beq
 J_Y|_{y(\ty)}\sseteq \tJ_Y,\quad  \quad
J_X|_{x(\tx,\ty)}\sseteq J_\tX+(\ty-\tf(\tx)),\quad \quad
  y(\ty)-f(x(\tx,\ty))\in (\ty-\tf(\tx))\sset \k\bl \tx,\ty\br.
\eeq
  Then  $y\to y(\ty)$ defines the needed  morphism of the targets, $\Phi_Y: \tY\to Y.$
  Substituting $\ty=\tf(\tx)$ we get $x\to x(\tx,\tf(\tx)),$
  which is the needed  morphism of the source, $\Phi_X :\tX\to X.$
    Recall that the ring $\k\bl x\br$ admits substitutions, \S\ref{Sec.Background.Schemes.Maps}.ii.
       Thus the morphisms on the diagram \eqref{Eq.diagram} are algebraic.

   Altogether we get:  $y(\ty)-f(x(\tx,\tf(\tx)))\in (\ty-\tf(\tx))\sset \k\bl \tx,\ty\br.$
 \item (\LRAPP) Apply the nested P\l oski approximation, \S\ref{Sec.Background.AP.classical}, to get a nested algebraic solution
 \beq\label{Eq.inside.proof.nested.Ploski}
 y(\ty,z_Y)\in \k\bl \ty,z_Y\br, \quad  x(\tx,\ty,z_X)\in \k\bl \tx,z_X\br, \quad \text{such that:}\quad
\eeq
\[ \hy(\ty)=y(\ty,\hz_Y(\ty))  \text{ and  }  \hx(\tx)=x(\tx,\ty,\hz_X(\tx,\ty))\ \text{for some }
 \hz_Y(\ty)\in \k[\![\ty]\!] \text{ and } \hz_X(\tx,\ty)\in \k[\![\tx,\ty]\!].
 \]
  Specialize the $x$-part, $x(\tx,\ty,z_X)\rightsquigarrow x(\tx,\tf(\tx),z_X).$
   Then $ y(\ty,z_Y)$ and $x(\tx,\tf(\tx),z_X)$ is the needed parametrization of the solution of \eqref{Eq.inside.proof.LRAP.Nash.Case}.
 \item (\LAP, \LAPP) Identify $R_X=R_\tX$ and impose the condition $\Phi_X=Id_X,$  cf. Remark \ref{Rem.LRAP}.iii.
  Equation  \eqref{Eq.inside.proof.LRAP.Nash.Case}    becomes: $\psi(J_Y)\sseteq J_\tY ,$ $\psi(y)-f(x)\sseteq J_X+(\ty-\tf(x)),$ with a nested solution $\hy(y).$ Proceed as before.
\item The inverse Artin approximation and its P\l oski version follow now by Remark \ref{Rem.LRAP}.iii.
 \epr
\eee

\subsection{Approximations for finite analytic morphisms,  $R_X\! =\!\quots{\k\{ x\}}{J_X},$   $R_Y\!=\!\quots{\k\{ y\}}{J_Y}$}\label{Sec.LRAP.finite.maps}

\bthe\label{Thm.LRAP.finite.analytic.maps}
 The properties \LRAP, \LRAPP, \LAP, \LAPP, the inverse Artin approximation (and its P\l oski version)  hold for   finite  analytic morphisms,   $ \Maps.$
\ethe
\noindent Namely, if $f,\tf\in \Maps$ are finite  analytic morphisms, and $\hPhi_Y\circ \tf=f\circ\hPhi_X,$ then $\hPhi_Y,\hPhi_X$ are approximated by  analytic morphisms (resp. are P\l oski-specializations, and so on).
\\\bpr
 Suppose the map-germs $f,\tf:X\to Y$ are formally $\cL\cR$-equivalent.

 We have   finitely-generated modules   $R_X\in mod$-$ R_Y$ (via $f$) and  $R_\tX\in mod$-$ R_\tY$ (via $\tf$).
   The (formal) relation $\hPhi_Y\circ \tf=f\circ\hPhi_X$ implies the (formal)   morphism of these
 modules with the base change, as on the diagram \eqref{Eq.diagram}.\vspace{-0.1cm}

\bee[\!\!\bf Step 1.]
\item
 Fix some generators of the modules, $R_X =R_Y\bl\{v_\bullet\}\br,$ and
 $R_\tX =R_\tY\bl\{\tv_\bullet\}\br. $
    W.l.o.g. we assume:   the set  $\{v_\bullet\}$ contains the generators $x=(x_1\dots x_n)$ of $R_X,$ while the set $\{\tv_\bullet\}$  contains the generators $\tx=(\tx_1\dots \tx_\tn)$ of $R_\tX.$
   For these sets of generators we get the presentation matrices of the modules,   $R_X=coker[A(y)]$, $R_\tX=coker[\tA(\ty)]$.

 The morphism of modules with base change, $R_Y \to R_\tY$ and $coker[A(y)]\to coker[\tA(\ty)],$    yields the commutative diagram
 \beq\bM
R^{\oplus (\dots)}_Y&\stackrel{A(y)}{\to}&R^{\oplus (\dots)}_Y&\to& R_X&\to& 0
\\
\downarrow U&&\downarrow V&& \quad \downarrow\Phi_X
\\R^{\oplus (\dots)}_\tY&\stackrel{\tA(\ty)}{\to}&R^{\oplus (\dots)}_\tY&\to& R_\tX&\to& 0
\eM \eeq
 Thus we get the condition
 $ V(\ty)\cdot A(y(\ty))= \tA(\ty)\cdot U(\ty)$, for the unknowns $y(\ty),$ $U(\ty),V(\ty)$
 over $R_\tY$.
 This is a system of implicit function equations, linear in $U,V,$   and analytic   in $y$.

The formal relation $\hPhi_Y\circ \tf=f\circ\hPhi_X$ yields a formal solution to this system, $\hy(\ty),\hU(\ty),\hV(\ty)$.
 Invoking the   Artin approximation over the analytic ring $R_\tY,$ \S\ref{Sec.Background.AP.classical}.i., we approximate this formal
 solution by an analytic solution, $y(\ty)$, $U(\ty)$, $V(\ty)$.
 This gives an analytic  morphism of modules with base change, $\Psi: coker[A(y)]\to coker[\tA(\ty)].$
 But this  morphism of  modules does not necessarily lift to an isomorphism of  algebras of the diagram \eqref{Eq.diagram}.

\item
 To ensure the lifting from modules to algebras we impose additional equations.
   Besides the relations of a module  we have the multiplications of the algebra:
   $\{\prod_{\sum\um=\ud}\uv^\um=\prod_{\sum\um'=\ud}\uv^{\um'}\}$.
    (Here $\ud\in \N^n$ is any multi-index.)
 We ask that the   morphism of modules $\Psi:coker[A(y)]\to coker \tA[\ty]$
 be multiplicative. This gives the equations:
 \beq
\Psi(\prod_{\sum\um=\ud} \uv^\um)=\prod_{\sum\um'=\ud} \Psi(\uv^{\um'})\in R_\tX,\quad \forall\ \ud\in \N^n
 \eeq
  These conditions are  polynomial equations in the entries of $V.$
 The set of these equations is infinite, but the ideal of these equations is generated by a finite subset.
  (Because $R_\tY$ is Noetherian, while $R_\tX$ is a finitely generated $R_\tY$-module.)
We add this (finite) set of equations to the equations $y(J_\tY)=J_\tY,$ $ V(\ty)\cdot A(y(\ty))= \tA(\ty)\cdot U(\ty)$.

A remark: the ideal $J_X\sset \k\{x\}$ represents zero in the $R_Y$-module $R_X.$ Hence it is sent to zero in $R_\tX.$


 Altogether we get a (finite) system of implicit-function equations.
     They are analytic power series in the unknowns $y,U,V$.
The formal morphism of the algebras $(R_X,R_Y,f)$, $(R_\tX,R_\tY,\tf)$ is a formal solution to these equations,
$\hU(\ty),\hV(\ty),\hy(\ty)$.
\bei
\item (\LRAP)
 By the Artin approximation, \S\ref{Sec.Background.AP.classical}.i.,  this solution is approximated by an analytic solution. And this ordinary solution is the needed analytic  morphism  of algebras.
 \item (\LRAPP)Using P\l oski theorem we present the formal solution as a specialization of a parameterized  analytic solution,
  $U(\ty,z),$ $V(\ty,z),$ $y(\ty,z)$.
\item (\LAP, \LAPP) Identify $R_X=R_\tX$ and impose the condition $\Phi_X=Id_X.$ This means: $Im(V-Id)\sseteq Im(\tA),$ i.e.
 $V-Id=\tA\cdot\tB$ for an unknown matrix $\tB.$
 We are still left with implicit function equations.
   \item  The inverse Artin approximation and its P\l oski version  follow by Remark \ref{Rem.LRAP}.
\epr
\eei
\eee

\section{Morphisms  of weakly-finite singularity type}

\subsection{Preparations}\label{Sec.Weak.Fin.Preparations} Below $R_X,R_Y$ are as in \S\ref{Sec.Background.Schemes.Maps}.
\bee[\!\!\!\!\!\!\bf i.\!]
\item The image of a morphism $X\stackrel{f}{\to}Y$ is the subscheme $f(X)\sseteq Y$ defined by the ideal $ker[f^\sharp_|: R_Y \to R_X]\sset R_Y.$
 A morphism is called dominant, if this ideal vanishes.
 In the analytic case one can also take small representatives of the germs $X,Y$ and to consider the set-theoretic image of $f.$
 For a finite morphism $f$ this set-theoretic image gives a (well-defined) analytic subgerm  coinciding with $f(X),$ e.g. by Grauert's theorem on direct image.

For non-finite morphisms this set-theoretic image can be pathological, e.g. not a well defined germ or not an analytic germ,   \cite{Joiţa-Tibar}.
Below we often refer to examples of Osgood-type: analytic morphisms $(\C^n,o)\to (\C^m,o),$   $2\le n<m,$ that are analytically dominant.
 (I.e. their set-theoretic image does not lie inside any analytic hypersurface, i.e. their image via $ker[f^\sharp]$ is the whole target.) See \cite[\S5.3]{Rond} and \cite{Abhyankar-van der Put.70} for examples of such pathologies.

 Take a small representative of such a morphism, $\C^n\supset Ball^n_\de\stackrel{f}{\to}Ball^m_\ep\sset\C^m.$ Then the germ of the    image can be different from  the  image of the germ, $ f(Ball^n_\de,x)\neq (f(Ball^n_\de),f(x))$. E.g. suppose $f$ is finite at a point $x\in Ball^n_\de,$ but of Osgood type at $o\in \C^n.$ Then $dim f(Ball^n_\de,x)=n,$ while
  $ (f(Ball^n_\de),f(x))=(\C^m,f(x)).$

\item
Take a (formal/analytic/Nash) germ $X$ and its module of  $\k$-linear derivations (vector fields), $Der_X:=Der_\k R_X.$
 E.g. for $X=(\k^n,o)$ one has $Der_{(\k^n,o)}:=R_{(\k^n,o)}\bl\frac{\di}{\di x_1},\dots,\frac{\di}{\di x_n}\br,$
 see Theorem 30.6(ii) of \cite{Matsumura}.
 For $X=V(J_X)\sset(\k^n,o)$ the elements of $Der_X$  lift to logarithmic derivations on $(\k^n,o)$ that preserve the ideal $J_X.$

More generally, given a subgerm $Z=V(I_Z)\sset X$ one considers the module of logarithmic derivations,
  $Der_{X,Z}:=Der_X Log(Z):=\{\xi|\ \xi(I_Z)\sseteq I_Z\}\sseteq Der_X.$
     The natural projection  $Der_{X,Z}\to Der_Z$ is not always surjective.
\\Yet more generally, for a set of subgerms $\{Z_j=V(I_{Z_j})\sset X\}_j$ we take ``multi-logarithmic" derivations:
\beq
Der_{X,\{Z_\bullet\}}:=\cap_j Der_X Log(Z_j)=\{\xi|\ \xi(I_{Z_j})\sseteq I_{Z_j},\ \ \forall j\}.
\eeq
 \item We use K\"ahler differentials. As our algebras are not f.g. over $\k,$ we work with separated differentials, $\Om^{sep}_{X/\k},$ see Exercise 16.14 of       \cite{Eisenbud-book}.
 E.g. $\Om^{sep}_{(\k^n,o)}=R_{(\k^n,o)}\bl dx_1,\dots,dx_n\br.$

 \item
 Suppose $\k$ is a local ring, then $f\in\Maps$ can be considered as a family of maps, over the base $Spec(\k).$
 Taking the maximal ideal, $\cm_\k\sset \k,$ one gets the central fibres:
  \bei
  \item of the schemes, $X_o=Spec(R_X\otimes \quots{\k}{\cm_\k})$ and $Y_o=Spec(R_Y\otimes \quots{\k}{\cm_\k}).$
 \item of the map,
 $f_o:=f\otimes \quots{\k}{\cm_\k}: X_o\to Y_o.$
  This $f_o$ is defined algebraically  via
  \beq
  f^\sharp_o:=f^\sharp\otimes \quots{\k}{\cm_\k}\in Hom_{\quots{\k}{\cm_\k}}(R_Y\otimes \quots{\k}{\cm_\k}, R_X\otimes \quots{\k}{\cm_\k}).
  \eeq
\eei
     If $\k$ is a field then (trivially) $f_o=f.$
 \eee

\subsubsection{Approximate integrability of derivations}\label{Sec.Weak.Fin.Integrability.of.Derivations} Take a (formal/analytic/Nash) germ $X\sseteq (\k^n,o)$ and a logarithmic derivation
  $\xi\in Der_{(\k^n,o),X}.$ Suppose $\xi(\cm)\!\sseteq\! \cm^2,$  i.e. this derivation is $\cm$-adically nilpotent. Define the self-map of $(\k^n,o)$ by $x\!\to\! x\!+\!\xi(x).$ This is an automorphism of $(\k^n,o),$ but it does not necessarily descend to an automorphism of $X,$ i.e. an automorphism of $R_X.$
  We want to get an automorphism of $X$ by correcting this map by ``higher order terms".
   As $\xi(\cm)\sseteq \cm^2,$ one has $\xi(\xi(\cm))\sseteq \cm\cdot \xi(\cm)\sseteq \cm^3.$
   \bed
   A logarithmic derivation $\xi\in  Der_{(\k^n,o),X}$ is called approximately integrable if
   there exists a (formal/analytic/Nash) automorphism $\Phi_{(\k^n,o)}\circlearrowright (\k^n,o)$ that descends to an automorphism $\Phi_X\circlearrowright X$ and satisfies:
    $\Phi_{(\k^n,o)}(x)-x-\xi(x)\in \xi(\xi(\cm)) .$
       \eed

More generally, for some subgerms $\{Z_\bullet\}\sset X$ take a multi-log-derivation $\xi\in Der_{ X,\{Z_\bullet\}}$ satisfying $\xi(\cm)\sseteq\cm^2$. Call $\xi$ {\em approximately integrable} if
  the automorphism      $\Phi_X\circlearrowright (\k^n,o)$ as above  defines automorphisms of $X$ and $\{Z_\bullet\}.$
\bex
\bee[\!\!\!\!\!\!\!\bf i.\!]
\item For $char(\k)=0$ every log-derivation (that satisfies $\xi(\cm)\sseteq \cm^2$) is  approximately integrable.
\\\bpr
Define $\Phi (x):=e^\xi x:=(\sum \frac{\xi^j}{j!})x,$ see \cite[Chapter II, \S6]{Bourbaki.Lie}. Thus $\Phi\circlearrowright \widehat{(\k^n\!,\!o)},$  a formal automorphism. It
 preserves  the subgerms
 \!$X,\!\{Z_\bullet\}$ and satisfies: \!$\Phi_X(x)\!-\!x\!-\!\xi(x)\!\in\! \xi(\xi(x)).$

 If needed, we approximate $\Phi $ by  an analytic/Nash automorphism that preserves  $X,\{Z_\bullet\}.$
  (The later conditions are implicit function equations, see Remark \ref{Rem.Background.Embedding.as.IF.eq}.)\epr

\item For $char>0$ the approximate  integrability is a non-trivial restriction, see e.g.
 \cite[\S2.5]{Kerner.Group.Orbits}.
\eee
\eex

\subsubsection{Approximate lifting of automorphisms}\label{Sec.Weak.Fin.Lifting.of.Automorphisms}
Take an ideal $I\sset R$ and an automorphism $\bPhi\in Aut_\k(\quots{R}{I}).$ (Below we assume $\bPhi(x)-x\in \cm^2,$ i.e. $\bPhi$ is ``tangent to identity".) It does not necessarily lift to an element of $Aut_\k(R).$
 But a weaker statement (used in the proof of Theorem \ref{Thm.LRAP.weakly.finite.analytic}) often holds:
\beq\label{Eq.Lifting.Property.Automorphisms}
\text{\em For any ideal $I\sset R$ and a suitable $N\gg1,$ any automorphism $\bPhi\in Aut_\k(\quots{R}{I^N})$}
\vspace{-0.1cm}
\eeq
\[
\text{\em has an approximate lifting, i.e. $\Phi\in Aut_\k(R)$ such that $\Phi(I)=I$ and $[\Phi]= [\bPhi]\in Aut_\k(\quots{R}{I}).$}
\]
(Here $N$ depends on $I\sset R$ only, not on $\Phi$.)
\\As was mentioned in \S\ref{Sec.Introduction.Results}, this property would follow from SAP for the filtration $I^\bullet\!\sset\! R.$ But this version of SAP does not hold, see
 \cite[Example 3.20.]{Rond}. For $char(\k)\!=\!0$ we give a different proof.
 \bel\label{Thm.Lifting.Property.Automorphisms}
  Property \eqref{Eq.Lifting.Property.Automorphisms} holds in zero characteristic.
 \eel
\bpr First we consider the formal case, $R=\quots{\k[\![x]\!]}{J}.$ The proof goes by ``iterations" on the (finite) number of generators of $J.$
 At each iteration we invoke the Artin-Rees lemma. Therefore we start with some $N\gg1,$ and by the end of iterations (Step 2) we specify conditions, which $N$ must satisfy. (These conditions depend on $I,J$ only, but not on $\bPhi.$)\vspace{-0.1cm}

\bee[\!\!\bf Step 1.]
\item
 Take the full preimage of $I$ in $\k[\![x]\!],$ denote it (abusing notations) by
 $I\sset \k[\![x]\!].$ Fix an automorphism $\bPhi\in Aut_\k(\quots{\k[\![x]\!]}{J+I^N}). $ As $char(\k)\!=\!0,$ we can present $\bPhi\!=\!exp(\bxi),$
  for a derivation $\bxi\!\in\! Der_\k(\k[\![x]\!])=:Der_{(\k^n,o)}$ satisfying: $\bxi(\cm)\sseteq \cm^2$ and $\bxi(J+I^N)\sseteq J+I^N.$  See \cite[Chapter II, \S6]{Bourbaki.Lie}.
\\We should approximate $\bxi$ by a derivation $\xi\!\in\! Der_{(\k^n,o)}$ satisfying: $ \xi(J )\!\sseteq\! J $ and $\xi\!\equiv\!\bxi\ \!mod\ \!I.$ (Then one defines $\Phi\!:=\!e^\xi.$) This approximation is constructed iteratively in Step 2.

\item
Fix some generators $J=(q_1,\dots,q_c)\sset \k[\![x]\!].$ Then $\bxi(q_1)\in J+I^N,$ i.e. $\bxi(q_1)-h_1\in I^N$ for some $h_1\in J.$
 Thus
 \beq
 \bxi(q_1)-h_1\in I^N\cap (Der_{(\k^n,o)}(q_1)+J)\stackrel{A.R.}{\sseteq }I^{N-d_1}\cdot (Der_{(\k^n,o)}(q_1)+J).
  \eeq
   Here   the inclusion $A.R.$ holds by the Artin-Rees lemma, for some $d_1.$ This $d_1$ is independent of $N,\bPhi,\xi.$
     Thus we can assume $N\gg d_1.$
\\
Therefore $\bxi(q_1)-\txi(q_1)\in J $ for some $\txi\!\in\! I^{N-d_1}\!\cdot\! Der_{(\k^n,o)}.$ Replace $\bxi$ by $\bxi-\txi$ to get:
\beq
\bxi(q_1)\in J, \quad  \text{  i.e. } \ \bxi\in Der_\k((q_1),J) , \quad \text{ and }\quad \bxi(q_2,\dots,q_c)\sseteq J+I^{N-d_1}.
\eeq
This adjusted derivation satisfies (as before): $\bxi(q_2)-h_2\in I^{N-d_1}\cap (Der_\k((q_1),J) \cdot q_2+J),$
 for some $h_2\in J.$ Invoking Artin-Rees we get: $\bxi(q_2)-h_2\in I^{N-d_1-d_2}\cdot (Der_\k((q_1),J) \cdot q_2+J).$
  Thus $\bxi(q_2)-\txi(q_2)\in J,$ for some $\txi\in I^{N-d_1-d_2}\cdot  Der_\k((q_1),J) .$ Replacing $\bxi$ by $\bxi-\txi$ we get:
  \beq
  \bxi(q_1,q_2)\in J, \quad  \text{  i.e. } \ \bxi\in Der_\k((q_1,q_2),J) , \quad \text{ and } \quad \bxi(q_3,\dots,q_c)\sseteq J+I^{N-d_1-d_2}.
  \eeq
 Iterate this up to the last generator $q_c.$ Eventually one gets a derivation $\xi\in  Der_{(\k^n,o)}$ satisfying $\xi(J)\sseteq J.$
 And this adjusted derivation $\xi$ differs from the initial one by an element of $I^{N-d_1-\cdots- d_c}\cdot Der_{(\k^n,o)}.$
\\\mbox{These integers $d_1\!\dots\! d_c$ depend on $I,\!J\!\sset\! \k[\![x]\!]$ only, not on $N,\!\xi,\!\bPhi.$ Thus it is enough to choose}
 $N\!>\!d_1 + \cdots + d_c,$ to get: \!$\forall\ \txi$ there exists $\xi\!\in\! Der_{(\k^n,o)}Log(J)$ such that $\txi\!-\!\xi\!\in\! I\!\cdot\! Der_{(\k^n,o)}.$
\item
Define the needed lifting of $\bPhi,$ i.e. the (formal) automorphism $\Phi:=e^\xi\in Aut_\k  ( \k[\![x]\!] ).$
 It satisfies the claimed conditions:
 \beq\label{Eq.inside.proof.Lifting.Autos}
  \Phi(J)=J, \quad \Phi(I)=I, \quad\quad \text{and }\quad  [\Phi]\equiv \bPhi\ mod\ I^{N-d_1-\cdots- d_c}.
  \eeq
\eee
In the analytic/Nash cases one applies the Artin approximation to the formal solution of  equation \eqref{Eq.inside.proof.Lifting.Autos}.
   (These  equations of are implicit function type.)
\epr

\bex\label{Ex.Lifting.Automorphisms}
Fix a chain of radical ideals, $I_1\sset \cdots\sset I_r\sset R,$ and the corresponding chain of projections, $R\to\quots{R }{I_1}\to\cdots\to \quots{R}{I_r}.$ Given an automorphism $\bPhi\circlearrowright \quots{R}{I_r}$ we want to lift it to an automorphism $\Phi\circlearrowright R$ that preserves the chain $I_\bullet$. Such a lifting does not exist in general (as before). However,   condition  \eqref{Eq.Lifting.Property.Automorphisms} implies:
\\{\em There exist  $N_1\dots N_r\gg1$ (depending on $R,I_\bullet$ only) such that any automorphism $\bPhi\!\circlearrowright\! \quots{R}{I^{N_1}_1\!+\!\cdots\!+\!I^{N_r}_r}$ lifts to an automorphism $\Phi\circlearrowright R$ that preserves the chain  $I_\bullet.$}

\bpr
 Fix (successively) numbers $1\ll N_1\ll N_2\ll\cdots\ll N_r,$ such that condition \eqref{Eq.Lifting.Property.Automorphisms}
  holds for the projections $R\to\quots{R}{I^{N_1}_1},$ $\quots{R}{I^{N_1}_1}\to \quots{R}{I^{N_1}_1+I^{N_2}_2},$ \dots, $\quots{R}{\sum^{r-1}_{j=1}I^{N_j}_j}\to \quots{R}{\sum^{r}_{j=1}I^{N_j}_j}.$ Fix any automorphism  $\Phi_r\circlearrowright\quots{R}{\sum^{r}_{j=1}I^{N_j}_j},$ lift it to $\Phi_{r-1} \circlearrowright \quots{R}{\sum^{r-1}_{j=1}I^{N_j}_j},$ and so on, up to $\Phi\circlearrowright R.$

    Observe: $\Phi_{r-1}(I_r)=I_r$ and $\quots{R}{I_r+\sum^r_{j=1}I^{N_j}_j}\circlearrowleft [\Phi_r]=[\Phi_{r-1}].$
     Similarly $\Phi_{r-2}(I_{r-1})=I_{r-1}$ and $\quots{R}{I_{r-1}+\sum^{r-1}_{j=1}I^{N_j}_j}\circlearrowleft [\Phi_{r-1}]=[\Phi_{r-2}].$
      Therefore $\Phi_{r-2}(I_r)=I_r.$ And so on\dots
\epr
\eex

\subsection{The critical locus/discriminant of  a   morphism $X\stackrel{f}{\to}Y$}\label{Sec.Weak.Fin.Critical.Locus}
Below we tacitly replace the target $Y$ by the image $f(X)\sseteq Y.$ Thus $f$ becomes a dominant morphism.
\bed\label{Def.Critical.locus.and.Discrim}
\bee
\item The critical module, $\cC_{X,f},$ is defined by the (exact) Zariski-Lipman sequence $Der_X f\to T_{X\to Y}\to \cC_{X,f}\to 0.$
 Here $T_{X\to Y}:=Hom_X(f^*\Om^{sep}_Y,R_X)$ is the tangent to the space $\Maps.$
\item The critical locus, $Crit_X f\sseteq X,$ is defined by the zeroth Fitting ideal $Fitt_0(\cC_{X,f})\sseteq R_X.$
\item The discriminant $\De\sset Y$ is the image $f(Crit_X f)\sseteq Y,$ see \S\ref{Sec.Weak.Fin.Preparations}.i.
\eee
\eed
For  questions of Artin approximation the particular scheme structures on $Crit_X f$ and $\De$ are not important. Usually we take the reduced structure.
 \bex\label{Ex.Critical.Locus} Take a dominant morphism $X\stackrel{f}{\to}Y.$
\bee[\!\!\!\!\!\bf i.\!]
\item  Let $Y=(\k^m,o),$ thus $\Om^{sep}_{Y}\cong R_Y\bl  dy_1,\dots,dy_m\br.$ As $f$ is dominant, $ker[f^\sharp:R_Y\to R_X]=0.$ Thus
 $f^*\Om^{sep}_{Y}$ is a free $R_X$-module.
Then $T_{X\to Y} \cong R^{\oplus m}_X$ and
  $\cC_{X,f}\cong\quots{R^{\oplus m}_X}{Der_X(f)}.$ Denote by $[Der_X(f)]$ the defining matrix of the submodule $Der_X(f) \sseteq R^{\oplus m}_X.$ Then $Crit_X f\sseteq X$ is defined by the determinantal ideal $I_m[Der_X(f)]\sseteq R_X.$
\\
For $X\!=\!(\k^n,o)$  this is the classical Jacobian matrix and $I_m[Der_{(\k^n,o)}(f)]$ is the classical Jacobian ideal.

 \item Let $Y\sseteq (\k^m,o)$ be defined by the ideal $I_Y\sset R_{(\k^m,o)}.$  Take the co-normal sequence
 \beq\label{Eq.CoNormal.Sequence}
     \quots{I_Y}{I^2_Y}\!\stackrel{d}{\to}\! R_Y\otimes \Om^{sep}_{(\k^m,o)}\!\to\! \Om^{sep}_Y\!\to\! 0.
\eeq
Then $Hom_X(f^*\Om^{sep}_Y,R_X)\!\sseteq\! Hom_X(f^*\Om^{sep}_{(\k^m,o)},R_X)\!\cong\! R^{\oplus m}_X.$

   Explicitly, fix some generators $I_Y=R_Y\bl \{q_\bullet\}\br$  to get:
   \beq\label{Eq.inside.ex.T_X/Y}
   T_{X\to Y}=\{h\in R^{\oplus m}_X|\ dq_1|_f(h)=0=\cdots=dq_c|_f(h)\}\sset R^{\oplus m}_X.
    \eeq
    Here $dq_i|_f \in f^\sharp\Om^{sep}_{(\k^m,o)}$ and $dq_i|_f(h)\in R_X.$

\item   Take a map $f\!:(\k^1,o)\!\to\! (\k^m,o),$ whose image is a curve-germ, $Y\!:=f(\k^1,o)\!=\!V(q_1,\dots,q_c).$
 Suppose $\k$ is infinite.
 If $char\!>\!0 $ we assume $(dq_1|_f,\dots,dq_c|_f)\!\neq\! (0,\dots,0).$ Then $\cC\!\cong\! \quots{R_{(\k^1,o)}}{(f'_1,\dots,f'_m)}.$

\bpr
 In this case:   $T_{ ( \k^1,o)\to Y}=\{h|\ dq_1|_f (h)=0=\cdots=dq_c|_f(h)\}\sset R^{\oplus m}_{(\k^n,o)}.$
   Take the generic coordinates in $(\k^m,o),$ so that
  $ord_x (\di_{y_1}q_j|_f)=\cdots  =ord_x( \di_{y_m}q_j|_f)$ for each  $j.$  Then for any $h_1\in R_{(\k^1,o)}$ the system
   $dq_1|_f(h)=0=\cdots=dq_c|_f(h)$ has a unique solution. Denote this solution by $(h_2,\dots,h_m).$
   Thus the projection $h=(h_1,\dots,h_m)\to h_1$ defines  the isomorphism $T_{  ( \k^1,o)\to Y}\isom{} R_{(\k^1,o)}.$
    It sends $Der_{(\k^1,o)}f\sset R^{\oplus m}_{(\k^1,o)}$ to $(f'_1)\sset R_{(\k^1,o)}.$
   By the genericity $ord(f'_1)=\cdots=ord(f'_m).$
   Hence the statement.
\epr

   \item (Osgood-type examples, $char\!=\!0$, cf. \S\ref{Sec.Weak.Fin.Preparations}) Suppose a non-finite analytic map $X\!\stackrel{f}{\to}\!(\k^m,o) $   is dominant, i.e. $
  \k\{y\}\!\supset\! ker[\k\{y\}\!\to\! R_X]\!=\!0,$ i.e. the set-theoretic image of (a small representative of) $f$ does not lie inside any analytic hypersurface.
  Suppose moreover $dim(X)\!<\!m.$ Thus $rank Der_X\!<\!m.$  Then the critical module $\cC_{X,f}\!=\!\quots{R^{\oplus m}_X}{Der_X(f)}$ is supported on the whole $X.$ Hence $Crit_X f\!=\!X.$

 \item When $char\!=\!p,$ the   pathology $Crit_X f\!=\!X$ can occur due to ramification. E.g. for a map $f:X\!\to\!(\k^m,o)$ of the form $f\!=\!(\tf^p_1,\dots,\tf^p_m),$ one has       $Crit_X f\!=\!X.$ Note that $f$ can be finite.

 \item Given a map $f:X\to  (\k^m,o),$ extend it to a map $\tf:=(f_1,\dots,f_m,\dots,f_{m+c}):X\to (\k^{m+c},o).$ Suppose $f_{m+1},\dots,f_{m+c}\in \k[f_1,\dots,f_m],$ then the projection $\tf(X)\to f(X)$ is an isomorphism. And thus $Crit_X f=Crit_X \tf .$

     \item A trivial case is the zero map,  $X\stackrel{f}{\to}\{o\}.$ Then $\cC_{X,f}=0$ and $Crit_{X,f}=X.$
\eee
\eex

\subsubsection{Basic properties  of $\cC_{X,f}$ and $Crit_X f$}\label{Sec.Weak.Fin.Critical.Locus.Properties}
\bee[\!\!\!\!\!\bf i.\!]
\item The ramification locus of a morphism $X\!\stackrel{f}{\to}\!Y$ is defined via the Fitting ideal of relative differentials, $Fitt_d[\Om^{sep}_{X/Y}].$ Here \!$d$=the dimension of the generic fibre. In the \!$\k$-smooth case, \!$X\!=\!(\k^n,o), Y\!=\!(\k^m,o),$  the ramification locus coincides with  the critical locus,
     $Fitt_0(\cC_{X,f})\!=\!Fitt_d[\Om^{sep}_{X/Y}],$ see Proposition 1.3 of \cite{Kallstrom}. \!Otherwise they differ.

   \item The critical module   $ \cC_{X, f} $ and  the critical   locus $Crit_X f $ are independent of the choice of coordinates in the target $(\k^m\!,o).$ Moreover, they do not depend on the choice of embedding $Y\!\hookrightarrow\! (\k^m,o).$
\\\bpr The minimal embedding of $Y$ is unique up to isomorphism, while any non-minimal embedding is equivalent to the extension of a minimal embedding by zeros.\epr

       \medskip

\noindent Any automorphism $\Phi_X\!\in \!Aut_X$   induces the equality: \mbox{$Crit_X(f\!\circ\!\Phi_X)\!=\!\Phi_X(Crit_X f)\!\sset\! X.$}

\item (Localization at a point) Let $\k=\bk,$ e.g. $\k=\C.$ Take a map of analytic germs, $f:X\to Y,$ and fix some representatives of these germs.
 Suppose $f(X,x)=(f(X),f(x))$ for each $x\in X,$  see \S\ref{Sec.Weak.Fin.Preparations}.i.
 Then $\cC_{X,f}$ becomes a coherent sheaf on  $X.$ For each (closed) point $x\in X$ we can take the analytic stalk $ (\cC_{X,f})_x.$
 Alternatively we can start from the analytic  map-germ $(X,x)\!\stackrel{f_x}{\to}\!(Y,f(x)),$ and get the critical module $\cC_{(X,x),f_x}.$
 \\{\em Claim:}   $\cC_{(X,x),f_x}\cong (\cC_{X,f})_x.$
\\\bpr The sheaves of  modules $Der_X f, T_{X\to Y}$ of definition \ref{Def.Critical.locus.and.Discrim}
   are functorial \wrt analytic localization.\epr

 \item ($Crit_{X}f$ measures the non-triviality of the fibration) The map of analytic spaces $f:X\to Y$ is called a locally trivial fibration at $x\in X$ if $(X,x)\cong (Y,f(x))\times B$ (isomorphism of analytic germs), with $f$ acting as the projection onto the first factor.
      Suppose $f(X,x)=(f(X),f(x)) $ for each $x\in X,$  see \S\ref{Sec.Weak.Fin.Preparations}.i.
      \\{\em Claim:}  Then
 $Crit_{X}f\sseteq \{x\in X|\ (X,x)\to (Y,f(x))\ \text{is not a locally trivial fibration}\}.$
\\\bpr
      If $(X,x)\to (Y,f(x))$ is   a locally trivial fibration, then
      \beq
      T_{ (X,x)\to(Y,f(x)) }\cong R_{(X,x)}\otimes Der_{(Y,f(x))} \text{  and }
       Der_{(X, x)}(f)\cong R_{(X,x)}\otimes  Der_{(Y,f(x))}(f).
       \eeq Thus $\cC_{(X,x),f}=0.$
      \epr
\eee

\subsubsection{A way to compute $Crit_X f$ via coverings}\label{Sec.Weak.Fin.How.to.Compute.Crit} Take a dominant morphism
 $f\!:\!X\! \to\!\! Y\!\sset\!\! (\k^m,\!o),$ and the projection $(\k^m,o)\!\to\!(\k^{m'},o) ,$
 $(y_1,\dots,y_m)\!\to \!(y_1,\dots,y_{m'}).$ Suppose it restricts to a (finite, surjective) covering $Y\stackrel{\pi}{\to} (\k^{m'},o).$
 We get the truncated dominant map $\pi\circ f: X\to (\k^{m'},o),$  whose critical locus is easy  to work with, by Example \ref{Ex.Critical.Locus}.i.
   Take the ramification locus, $Ram(\pi)\sset Y,$ see \S\ref{Sec.Weak.Fin.Critical.Locus.Properties}.i., and its preimage  $f^{-1}(Ram(\pi))\sset X.$
\bel\label{Thm.Crit.computing.via.finite.covering}
Suppose $f^{-1}(Ram(\pi))\sset X$ contains no irreducible components of $Crit_X(f) $ or of $Crit(\pi\circ f).$ Then
 $Crit_X f=Crit_X(\pi\circ f)\sset X.$
\eel
\bpr
We compare the critical modules $\cC_{X,f}$ and $\cC_{X,\pi\circ f }.$ The truncation $\pi$ induces the diagram
\beq\bM
 Der_X f&\to &T_{X\to Y}&\to& \cC_{X,f}&\to& 0
\\\downarrow \pi_1& & \downarrow \pi_2&& \downarrow \pi_3
\\Der_X (\pi\circ f)&\to& T_{X\to (\k^{m'},o)}&\to& \cC_{X,\pi\circ f}&\to& 0
\eM
\eeq
Here the maps $\pi_1,\pi_2$ define $\pi_3,$ a homomorphism of $R_X$-modules.

We claim: $\pi_1,$ $\pi_2$ are  isomorphisms over $\cU:=X\smin f^{-1}(Ram(\pi)).$
Indeed fix some generators $I_Y=R_Y\bl q_1,\dots,q_c \br,$ thus $q_\bullet(f)=0.$
 Any derivation $\xi\in Der_X$ satisfies: $dq_\bullet|_f(\xi(f))=\xi(q_\bullet(f))=0.$
  Over $\cU$ the matrix $[\di_{m'+1}q_\bullet|_f,\dots,\di_m q_\bullet|_f]\in Mat_{c\times(m-m')}(R_X)$ is of   constant rank $m-m'.$
 Therefore  over $\cU $ we get:
 \bei
 \item $Der_\cU(f_{m'+1}),\dots,Der_\cU(f_m)\sseteq Der_\cU(f_1,\dots,f_{m'}) .$ Hence $\pi_1$ is an isomorphism.
 \item For any $(h_1,\dots,h_{m'})\in R^{\oplus m'}_\cU$ there exist  (and unique) $h_{m'+1},\dots,h_{m }\in R_\cU$ satisfying:
  $dq_i|_f(h)=0$ for any $i.$ Hence $\pi_2$ is an isomorphism, see \eqref{Eq.inside.ex.T_X/Y}.
 \eei
Therefore $\pi_3$ is an isomorphism over $\cU .$

 Therefore  $Crit_X f\cap\cU=Crit_X(\pi\circ f)\cap\cU\sset X.$ Finally:
 \beq
Crit_X f\stackrel{!}{=}\overline{ Crit_X f\cap\cU}=\overline{Crit_X(\pi\circ f)\cap\cU}
\stackrel{!}{=}Crit_X(\pi\circ f)\sset X.
 \eeq
The equalities (!) hold by the initial assumption.
\epr

 \subsubsection{Restricting to the central fibre}\label{Sec.Weak.Fin.Restricting.to.Central.Fibre}
        For $(\k,\cm_\k)$ a local ring, and $X\stackrel{f}{\to}Y$ a dominant morphism,
        we  take the central fibre $X_o\stackrel{f_o}{\to}Y_o$, see \S\ref{Sec.Weak.Fin.Preparations}.iv.
        One can either compute the critical locus of the central fibre, $Crit_{X_o}f_o,$ or take the central fibre of the critical locus, $Crit_X f\cap V(\cm_\k).$
 Under certain assumptions they coincide. We prove a more general statement.

Let $\k$ be a local Cohen-Macaulay algebra (over a field $\k_o$), with a system of parameters $t=(t_1,\dots,t_c).$
Suppose this sequence is also $R_X$-regular.
Denote $X_o\!=\!Spec(\quots{R_X}{(t)}),$ $Y_o\!=\!Spec(\quots{R_Y}{(t)}).$ Thus $f(X_o)\sseteq f(X)\cap Y_o=Y_o,$ the images in the sense of  \S\ref{Sec.Weak.Fin.Preparations}.i.
  But the inequality can occur.
\bex
Let $\k=\k_o[\![t]\!]$ and $X=X_o\times(\k^1_o,o),$ $Y=(\k^m_o,o)\times(\k^1_o,o).$ Fix a dominant morphism $X_o\stackrel{g}{\to}(\k^m_o,o),$
 and define $f:X\to Y$ by $(x_o,t)\to(t\cdot g(x_o),t).$ Then $f$ is dominant as well. But $f_o(X_o)=\{o\}\times \{o\}\ssetneq f(X)\cap Y_o=Y_o .$
\eex
\bel
 If $X_o\stackrel{f_o}{\to}Y_o$  is dominant, 
 then $\cC_{X_o,f_o}=\quots{\k}{(t)}\otimes\cC_{X,f}$ and
$Crit_{X_o}f_o=Crit_X f\cap X_o.$
\eel
 \bpr
To obtain the critical locus $\cC_{X_o,f_o}$ we should restrict $Der_X$ to $Der_{X_o}$ and $T_{X\to Y} $
 to $T_{X_o\to Y_o}. $ In both cases we use an auxiliary fact, recalled in Step 1.\vspace{-0.1cm}

 \bee[\!\!\bf Step 1.]
\item  Take a homomorphism of rings $R\to S,$   a f.g. module $M\in mod$-$R,$ and its Auslander transpose $M^t.$
 In Step 2 we use the (standard) exact sequence:
 \beq\label{Eq.Auslander.dual}
 0\to Tor^R_2(M^t,S)\to S\otimes Hom_R(M,R)\to Hom_S(S\otimes M,S)\to0.
 \eeq
For completeness we recall its proof. Fix a presentation of the module, $F_1\stackrel{A}{\to}F_0\to M\to0.$ Apply $Hom_R(-,R)$ to get
 $0\to M^*\to F^*_0\stackrel{A^t}{\to}F^*_1\to M^t\to 0.$ Split this sequence into two
 \beq
 0\to M^*\to F^*_0\stackrel{A^t}{\to} Im(A^t)\to0 \quad\quad\quad \text{and}\quad\quad \quad 0\to Im(A^t)\to F^*_1\to M^t\to 0.
 \eeq
Apply $S\otimes $ to the first sequence to get: $0\to Tor^R_1(Im(A^t),S)\to S\otimes M^*\to S\otimes F^*_0\stackrel{S\otimes A^t}{\to}S\otimes Im(A^t)\to0.$ Here $ker[S\otimes A^t]\cong(S\otimes M)^*.$
\\Apply $S\otimes $ to the second sequence to get: $0\to Tor^R_2(M^t,S)\to Tor^R_1(Im(A^t),S)\to 0.$
Combine this to get: $0\to Tor^R_2(M^t,S)\to S\otimes M^*\to  (S\otimes M)^*\to0,$ as claimed.
 \item
 It is enough to prove the claim for the restriction $R\to \quots{R}{(\tau)},$ where $\tau\in \cm_\k\sset \k$ is any regular element.
  Then  the claim for the restriction $R\to \quots{R}{(t) },$ with the system of parameters $t=(t_1,\dots,t_c),$  follows by induction.
 \\Thus we apply \eqref{Eq.Auslander.dual} for $S\!:=R_{X_o}\!:=\quots{R}{(\tau)}.$ In this case  $0\!\to\! R_X \! \stackrel{\tau}{\to}\!R_X\!\to\! R_{X_o}\!\to\!0$ is exact.
    Therefore   $Tor^R_2(-, R_{X_o})\!=\!0.$
     Hence $R_{X_o}\!\otimes  Hom_R(M,R) \isom{}   Hom_{R_{X_o}}(R_{X_o} \otimes  M,R_{X_o}).$
\bei
\item
For   $M=\Om^{sep}_X $ one gets: $  R_{X_o}\otimes Der_X\isom{} Der_{X_o} .$

\item
For  $M=f^*\Om^{sep}_Y $ one gets: $  R_{X_o}\otimes T_{X\to Y}\isom{} T_{X_o\to Y_o} .$
\eei
Therefore, applying $ R_{X_o}\otimes$ to $Der_X\to T_{X\to Y}\to \cC_{X/Y}\to0,$ one gets:
$Der_{X_o}\to T_{X_o\to Y_o}\to R_{X_o}\otimes\cC_{X,f}\to0.$  (Recall that both $f$ and $f_o$ are dominant.)

Hence  $\cC_{X_o,f_o}=\quots{R_X}{(\tau)}\otimes\cC_{X,f}.$ And thus $Crit_{X_o}f_o=Crit_X f\cap V(\tau).$
 \epr
\eee

\bex
Let $\k$ be a regular local algebra, and $R_X$ flat over $\k.$ Suppose $f(X_o)=f(X)\cap Y_o.$ Then
  $\cC_{X_o,f_o}=\quots{R_X}{\cm_\k\cdot R_X}\otimes\cC_{X,f}$ and
$Crit_{X_o}f_o=Crit_X f\cap V(\cm_\k).$
\eex

\subsection{Every dominant map  $X\!\stackrel{f}{\to}\!Y$ is finitely $\cR$-determined by its critical locus $\pmb{Crit_X f}$}
 For $char>0$ we assume: there exists $N\!\gg\!1,$ such that all the derivations in $Fitt_0(\cC_{X,f})^N\!\cdot\! Der_X,$     are approximately integrable, see \S\ref{Sec.Weak.Fin.Integrability.of.Derivations}.
\bel\label{Thm.Critical.Locus.Fin.Determ.Map}
  \mbox{There exists $N\!\gg\!1$ satisfying: any map $\tf\!:\! X\!\to\! Y$ for which $\tf\!-\!f\!\in\! Fitt_0(\cC)^N\!\cdot\! R_X^{\oplus m}$} is $\cR$-equivalent to $f.$
  (I.e. $\tf=f\circ\Phi_X$ for an automorphism $\Phi_X\in Aut_X.$) Moreover, one can choose $\Phi_X$ satisfying:  $\Phi_X(x)-x\in Fitt_0(\cC)^{N-d},$ for some $d$ that depends on $f$ only, but not $N,\tf.$
\eel
\bpr In the proof we use the Artin-Rees lemma twice, and then apply the integrability assumption. Therefore we fix $N\gg1,$ and give the precise bound at the end of Step 1.\vspace{-0.1cm}

\bee[\!\!\bf Step 1.]
\item Fix some generators of the defining ideal of the image, $I_Y=R_Y\bl\{q_\bullet\}\br.$ (Thus $q_\bullet(f)=0.$) Then
\beq\label{Eq.inside.proof.finite.R.determinacy}
0=q_i(\tf)=q_i(f+\tf-f)=dq_i|_f(\tf-f )+(\dots), \quad \text{ where }\quad (\dots)\in (\tf-f)^2.
\eeq
Take the preimage of ideal $dI_Y|_f=f^{-1}(dI_Y)\sset f^{-1}(R_Y\otimes \Om^{sep}_{(\k^m,o)})\sseteq R^{\oplus m}_X $
 (cf. Example \ref{Ex.Critical.Locus}.ii.)
 By \eqref{Eq.inside.proof.finite.R.determinacy} one has: $dI_Y|_f(\tf-f)\sset (\tf-f)^2\sseteq Fitt_0(\cC)^{2N}.$

 Moreover, take the matrix $[dI_Y|_f]=[dq_1|_f,\dots,dq_c|_f]^T\in Mat_{c\times m}$ and the corresponding
 map $ [dI_Y|_f]:  R^{\oplus m}_X\to R^{\oplus c}_X.$ The image of this map is a (finitely generated) submodule of $R^{\oplus c}_X.$ Therefore:
 \beq
 [dI_Y|_f](\tf-f)\in Fitt_0(\cC)^{2N}\cdot R^{\oplus c}_X\ \cap\ [dI_Y|_f](R^{\oplus m}_X)\stackrel{A.R.}{\sseteq} Fitt_0(\cC)^{2N-d'}\cdot[dI_Y|_f](R^{\oplus m}_X).
 \eeq
 The inclusion holds by the Artin-Rees lemma, and the integer $d'$ depends on $f$ only, but not on $\tf$ or $N.$ (In particular, $d'\ll N.$)

Thus $ [dI_Y|_f](\tf-f)= [dI_Y|_f](h)\sset R^{\oplus c}_X$ for some vector $h\in Fitt_0(\cC)^{2N-d'}\cdot  R^{\oplus m}_X .$
 Hence $ [dI_Y|_f](\tf-f-h)=0\in R^{\oplus c}_X.$ Therefore, by \eqref{Eq.inside.ex.T_X/Y}, $\tf-f-h\in T_{X\to Y}\cap  Fitt_0(\cC)^N\cdot  R^{\oplus m}_X.$ Applying Artin-Rees lemma to the (finitely generated) submodule $ T_{X\to Y}\sseteq  R^{\oplus m}_X,$ one has:
  $\tf-f-h\in Fitt_0(\cC)^{N-d+1}\cdot  T_{X\to Y}.$ Here the integer $d$  does not depend on $N,\tf.$ (Hence $d\ll N.$)
  And therefore  $\tf-f-h\in Fitt_0(\cC)^{N-d}\cdot  Der_X f.$
 Altogether:
 \beq\label{Eq.inside.proof.finite.R.determinacy.1}
 \tf=f+\xi(f)+h, \quad \text{ where }\xi\in Fitt_0(\cC)^{N-d}\cdot  Der_X, \quad   \text{ while }   h\in Fitt_0(\cC)^{2N-d'}\cdot R^{\oplus m}_X.
 \eeq
 A remark: the integers $d,d'$ are preserved under $Aut_X.$

 We assume $N\gg1,$ such that:
 \beq\label{Eq.condition.on.N}
  d'<N \  \text{ and } \ d<N, \text{ and all the derivations in $Fitt_0(\cC)^{N-d}\cdot Der_X$ are integrable.}
\eeq
\item
  By our  assumption $\xi$ can be approximately integrated to an automorphism $\Phi_X\in Aut_X$ satisfying: $\Phi_X(x)=x+\xi(x)+(\dots)$ with $(\dots)\in Fitt_0(\cC)^{N+1}.$
   Then $ \tf\circ \Phi^{-1}_X- f\in Fitt_0(\cC)^{N+1}\cdot R^{\oplus m}_X. $

Now we repeat Step 1 for $\tf\circ\Phi^{-1}_X.$ Note that $d',d$ are preserved. Thus the condition \eqref{Eq.condition.on.N}  remains true.

 Iterate this process to get: $\tf\circ \Phi^{-1}_N\circ\cdots\circ \Phi^{-1}_{N+j}\in \{f\}+Fitt_0(\cC)^{N+j}\cdot R^{\oplus m}_X$ for each $j\in \N.$
  Define the formal automorphism $\hPhi_X:=\lim (\cdots\circ\Phi_{N+1}\circ \Phi_N).$ This limit exists because
   $\Phi_{N+j}(x)-x\in Fitt_0(\cC_{X,f})^{N+j-d}$ for every $j\in \N.$

  We get: $\tf=f\circ \hPhi_X,$ i.e. $\tf,f$ are formally $\cR$-equivalent.
  And by our construction $\hPhi_X(x)-x\in Fitt_0(\cC)^{N-d}\cdot \hR_X.$ This proves the statement in the formal case.
\eee
  In the analytic/Nash cases we apply the Artin approximation, Example \ref{Ex.RAP}.
\epr

\bex
Let $\k$ be a field and suppose $X\stackrel{f}{\to}(\k^m,o)$ is dominant. Then $f$ is finitely determined by $Crit_{X}f\sset X.$
 This is classically known, see  \cite[\S4.1]{Kerner.Group.Orbits} for the references.
\eex

\subsection{Morphisms of   finite singularity type}\label{Sec.Weak.Fin.Finite.Sing.Type}
\bed\label{Def.Finite.Sing.Type}
A dominant morphism $f\in \Maps$ is called of finite singularity type if the restriction $f_|:Crit_X f\to \De$ is a finite morphism.
\eed
\bex
\bee[\!\!\!\!\bf i.\!]
\item
For $Y=(\k^m,o)$ this gives the classical notion. In particular:
\\
The restricted map $Crit_X f\stackrel{f|}{\to}\De$ is  finite   iff $f:X\to (\k^m,o)$ is contact ($\cK$)-finite.
  (And then, for a field $\k=\bk,$ the  germ $V(f)\sset X$ is either of zero dimension  or  an ICIS.)
 The $\C$-analytic case is well known,
   \cite[pg.224]{Mond-Nuno}, \cite{Mond-Montaldi}.
 For the case of arbitrary characteristic see   \cite[Lemma 3.15]{Kerner.Group.Orbits}.

\item This finiteness condition does not depend on the choice of
  scheme structure on $Crit_X f$, e.g. it is preserved if one passes to $Spec(\quots{R_X}{\sqrt{Fitt_0(\cC)}})$ or to
  $Spec(\quots{R_X}{  Fitt_0(\cC) ^N})$.

\item Under the assumptions of Lemma \ref{Thm.Crit.computing.via.finite.covering}: $X\stackrel{f}{\to}Y$ is of finite singularity type iff
 $X\stackrel{\pi\circ f}{\to}(\k^{m'},o)$ is such.

\item
 Suppose $Crit_X f=X,$ see Example \ref{Ex.Critical.Locus}. Then $f$ is  of finite singularity type iff it is a finite morphism.

\item (Preservation of f.s.t. in deformations) Let $\k$ be a Cohen-Macaulay local algebra over a field $\k_o.$ Suppose $R_X$ is flat over $\k.$
 Suppose the central fibre of a map, $X_o\stackrel{f_o}{\to}Y_o,$ is of finite singularity type.
Then the map $X \stackrel{f }{\to}Y $ is of finite singularity type.
\\\bpr
 The map $Crit_{X_o} f_o \stackrel{f_o|}{\to}\De_o$ is finite. It is dominant, thus $f$ is dominant.
 Hence  (\S\ref{Sec.Weak.Fin.Restricting.to.Central.Fibre})
     the map  $Crit_X f|_{X_o} \stackrel{f|}{\to}\De_o$ is finite. And therefore  the map  $Crit_X f \stackrel{f|}{\to}\De$ is finite as well.
    This follows by the Weierstrass finiteness theorem, which holds for W-systems, \cite[\S1]{Denef-Lipshitz}, in particular for our rings.
\epr

    \eee
\eex

\subsection{Higher critical loci}\label{Sec.Weak.Fin.Higher.Critical.Loci}

 Take a dominant morphism $f:X\to Y,$  its critical locus $Crit_1:=Crit_X f\sseteq X,$ and the discriminant $\De_1:=f(Crit_1)\sset Y.$
  (With reduced scheme structures.)
  One could restrict $f$ to this critical locus, $Crit_1\stackrel{f_|}{\to}\De_1,$ and study the critical locus of the restriction, $Crit_{Crit_1}(f_|).$ But this ignores the source and looses information about $X.$ E.g.,  for $X\neq (\k^n,o),$ not all derivations in $Der_{Crit}$ come from those of $Der_X.$

For our approximation results we need to keep $X.$ Hence the higher critical loci are defined (iteratively) via logarithmic drivations.
 Put $\cC_0:=R_X,$ $Crit_0:=X,$ $\De_0:=Y.$ Suppose we have  defined the
  triples $(\cC_0,Crit_0,\De_0),\dots, (\cC_j,Crit_j,\De_j).$  Note that all the restrictions $Crit_i\stackrel{f_|}{\to}\De_i$ are dominant.
 \bed
\bee
\item
 The $(j+1)$-critical module is defined by the sequence of $R_X$-modules
 \beq\label{Eq.def.higher.crit.locus}
R_{Crit_j}\otimes Der_{Crit_0\dots Crit_j}f\to T_{Crit_j\to \De_j}\to \cC_{j+1}\to0,
 \eeq
  where $Der_{Crit_0\dots Crit_j}=\cap^j_{i=0} Der_X Log (Fitt_0(\cC_i))$ and $T_{Crit_j\to \De_j}=Hom_{Crit_j}(f^*\Om^{sep}_{\De_j},R_{Crit_j}).$

\item The $(j+1)$-critical locus is the support of $\cC_{j+1}.$
 Namely, $Crit_{j+1}=Supp(\cC_{j+1})\sseteq X$ is defined by the Fitting ideal $Fitt_0(\cC_{j+1})\sset R_X.$
\item The $(j+1)$-discriminant is the image, $\De_{j+1}=f(Crit_{j+1})\sseteq Y,$ see \ref{Sec.Weak.Fin.Preparations}.i.
\eee
 \eed
 We get  the   higher critical loci and higher discriminants (the formal/analytic/Nash germs):
\beq
X=Crit_0\supseteq Crit_1 \supseteq  Crit_2 \supseteq\cdots,
\hspace{2cm}
Y=\De_0\supseteq \De_1 \supseteq \De_2 \supseteq\cdots.
\eeq
They are all taken with reduced structure, as in definition \ref{Def.Critical.locus.and.Discrim}.

\subsubsection{Basic properties of $Crit_\bullet$}\label{Sec.Weak.Fin.Higher.Crit.Loc.Properties}
 \bee[\!\!\!\!\bf i.\!]
 \item
 The subgerms $Crit_\bullet(f)\sseteq X$ depend neither on the choice of   embedding $(Y,o)\sseteq (\k^m,o)$ nor on the choice of coordinates in $(\k^m,o).$

  For any automorphism $\Phi_X\circlearrowright X$ one has  $Crit_i(f\circ\Phi_X)=\Phi_X(Crit_i(f))\sseteq X.$



\item Similarly, the discriminants $\De_\bullet\sset Y$ are preserved under the action of $Aut_X,$ while any automorphism $\Phi_Y\in Aut_Y$ satisfies:
 $\Phi_Y(\De_\bullet(f))=\De_\bullet(\Phi_Y(f)).$

\item    Suppose $Y\!=\!(\k^m,o)$ and $\De$ is smooth over $\k.$ Applying $Aut_Y$
 we rectify $\De\!=\!V(y_{ c+1},\dots,y_m)\!\sset \!(\k^m,o).$ Then we have the truncated morphism $f_|\!:=\!(f_1,\dots,f_c):Crit_X\!\to\! (\k^c,o)$ and $Crit_2$ contains the critical locus of $f_|.$ Moreover, if $X=(\k^n,o),$ then $Crit_2\!=\!Crit(f_|).$

\item More generally, for $f:X\to Y,$ fix some $j\ge1$ and consider the restriction $f_|:Crit_j\to \De_j.$ Then $\cC_{j+1,f}\twoheadrightarrow \cC_{1,f_|}.$ Therefore $Crit_{j+1}(f)\supseteq Crit_{1}(f_|).$

  Indeed, $R_{Crit_j}\otimes Der_{Crit_0\dots Crit_j}f\sseteq Der_{Crit_j}(f_|).$ Now invoke \eqref{Eq.def.higher.crit.locus}.

\item Take a map $f:(\k^n,o)\to Y.$ Then $Der_{(\k^n,o),Crit}\twoheadrightarrow Der_{ Crit}.$ Therefore for the restriction $f_|:Crit\to \De$ one has: $Crit_{2,f}=Crit_{1,f_|}.$

\item ($\k$ an infinite field) Take a dominant morphism $(\k^n,o)\!\stackrel{f}{\to}\!(\k^m,o),$ with $m\ge2.$ Suppose the central fibre $ V(f)\!\sset\! (\k^n,o)$ contains no irreducible component of $Crit(f).$
  Suppose $\De\!\sset \!(\k^m,o)$ is a hypersurface germ. Take   generic coordinates on $(\k^m,o).$ Then $Crit_2(f)\!=\!Crit_{Crit_1}(f_1\dots f_{m-1}).$

\bpr Take the projection $(\k^m,o)\supset\De\stackrel{\pi}{\to}(\k^{m-1},o),$ and its ramification locus $Ram(\pi).$ By our assumptions $f^{-1}(Ram(\pi))$ contains no irreducible components of $Crit(f), Crit(\pi\circ f).$ Thus $Crit_{Crit}(f_|)=Crit_{Crit}(f_1\dots f_{m-1})_|,$
 by Lemma \ref{Thm.Crit.computing.via.finite.covering}. Finally, $Crit_2=Crit_{Crit}(f_|),$ as $(\k^m,o)$ is $\k$-smooth.
\epr

\item Suppose $\De_j=\{o\}.$ Then $Crit_j=Crit_{j+1}=\cdots,$ see Example  \ref{Ex.Critical.Locus}.vii.
\eee

\subsubsection{Every dominant map is finitely $\cR$-determined by its higher critical locus} Fix some $j\ge0$ and suppose the  dominant maps $X\stackrel{f,\tf}{\to}Y$ satisfy:
 $Crit_i(f)\!=\!Crit_i(\tf)$ and $\De_i(f)\!=\!\De_i(\tf)$ for $i=0,\dots,j.$ (Here $Crit_0=X$.)
  For $char>0$ we assume: all the derivations in the modules
\beq\label{Eq.integrability.conditions}
I^N_{Crit_1} \cdot Der_{Crit_0}, \quad\quad\quad I^N_{Crit_2} \cdot Der_{Crit_0,Crit_1},\quad\quad\quad\dots,  \quad\quad\quad
 I^N_{Crit_{j+1}} \cdot Der_{Crit_0,\dots, Crit_j},
 \eeq
 are approximately integrable for $N\gg1,$ see \S\ref{Sec.Weak.Fin.Integrability.of.Derivations}.
\bel\label{Thm.Higher.Critical.Locus.Fin.Determ.Map}
 If $f\!-\!\tf\!\in\! \ I^N_{Crit_{j+1}}\!\cdot\! \RmX$ for some $N\!\gg\!1,$ then $\tf$ is $\cR$-equivalent to $f,$ i.e. $\tf\!=\!f\!\circ\!\Phi_X$ for $\Phi_X\!\in \!Aut_X.$ Moreover, one can ensure $\Phi_X(x)\!-\!x\!\in \! I_{Crit_{j+1}}^{N-d}.$

   Here $N$ depends on $f$ only (not on $\tf$),  and $d$  depends on $f$ only (not on $\tf,N$).
\eel
\bpr \mbox{We assume $N\!\gg\!1,$ and give the explicit bound at the end.
Repeat the proof of Lemma \ref{Thm.Critical.Locus.Fin.Determ.Map}}
 \mbox{for the restrictions  $Crit_j\!\stackrel{f_|,\tf_|}{\to}\!\De_j.$
 (Note that $Crit_{j+1}(f)\!\supseteq\! Crit_1(f_|),$ see \S\ref{Sec.Weak.Fin.Higher.Crit.Loc.Properties}.iv.)
 Equation \eqref{Eq.inside.proof.finite.R.determinacy.1} is:}
 \beq
 \tf_|=f_|+\xi(f_|)+h,\quad \text{for}\quad \xi\in Fitt_0(\cC_{j+1})^{N-d}\cdot Der_{Crit_0\dots Crit_j},\quad h\in Fitt_0(\cC_{j+1})^{2N-d'}\cdot \RmX.
 \eeq
Approximately integrate the derivation $\xi$ to an automorphism $\Phi_X\in Aut_X$ that preserves the chain
 $X\supset  Crit_1\supset \cdots\supset Crit_j$
 and satisfies: $\Phi_X(x)-x\in Fitt_0(\cC_{j+1})^{N-d}.$
  Then (as in Lemma \ref{Thm.Critical.Locus.Fin.Determ.Map}) $\tf_|\circ\Phi_X^{-1}-f_|\in Fitt_0(\cC_{j+1})^{N+1}\cdot \RmX. $

Iterate this (as in Lemma \ref{Thm.Critical.Locus.Fin.Determ.Map}), and apply \RAP (if needed) to get the automorphism $\Phi_X\in Aut_X$ that preserves the chain
 $X\supset  Crit_1\supset \cdots\supset Crit_j$ and satisfies: $\tf_|\circ\Phi_X^{-1}=f_|.$

 Replace $\tf$ by $\tf \circ\Phi_X^{-1}.$ Then the restrictions of $f,\tf$ to $Crit_j$ coincide.
  We get  $f -\tf   \in I_{Crit_j}\cdot \RmX.$

 Now repeat the whole process, restricting to   $N_j$-th infinitesimal neighborhood of $Crit_j\sset X.$ We get:
  $f -\tf   \in I_{Crit_j}^{N_j}\cdot \RmX.$ Here $N_j\gg1,$ accordingly one takes $N\gg N_j.$
\\
Iterate this argument for $Crit_{j-1} $ and up to $ Crit_0\!=\!X.$ At this step one has: $  f\! -\!\tf \!\in \! Fitt_0(\cC_0)\!\cdot\! \RmX\!=\!0.$

Finally, the lower bound on $N$ is: conditions \eqref{Eq.integrability.conditions} hold, and moreover, Lemma \ref{Thm.Critical.Locus.Fin.Determ.Map} is applicable for each restriction [$N_i$-th neighborhood of $Crit_i$]$\stackrel{f_|}{\to}$[$N_i$-th neighborhood of $\De_i$]. Thus $1\ll N_1\ll\cdots\ll N_j\ll N.$
\epr

\subsection{Morphisms of weakly-finite singularity type}\label{Sec.Weak.Fin.Weakly.Finite.Sing.Type}
  \bed\label{Def.Weakly.Finite.Sing.Type}
A (dominant) morphism $f\in \Maps$ is called of   \underline{weakly-finite singularity type} if
 the restriction $f_|:Crit_r f\to \De_r$ is a finite morphism for some $r\ge0$.
\eed
\bex
\bee[\!\!\!\bf i.\!]
\item
For $r=0$ this means: $f$ is a finite morphism.

 For $r=1$ we get just the finite singularity type, definition \ref{Def.Finite.Sing.Type}.

\item Being of weakly-finite singularity type does not depend on the choice of the embeddings $X\sseteq (\k^n,o),$ $Y\sseteq(\k^m,o),$ and is preserved under the group action $Aut_X\times Aut_Y\circlearrowright\Maps,$ see \S\ref{Sec.Weak.Fin.Higher.Crit.Loc.Properties}.

\item Let $f:X\to (\k^1,o) $ and suppose $\De=(o)\sset (\k^1,o)$. Then $f$ is of weakly-finite singularity type iff it is of finite singularity type. (See \S\ref{Sec.Weak.Fin.Higher.Crit.Loc.Properties}.vii.)

   \item
   Take a dominant morphism $f:(\k^n,o)\to (\k^m,o),$ suppose the (reduced) discriminant $\De\sset (\k^m,o)$ is $\k$-smooth. Rectify it to $\De=V(y_{c+1},\dots,y_m)\sset (\k^m,o).$
   If the truncation $(f_1,\dots,f_c):\ Crit_{X}f\to (\k^c,o)$ is of finite singularity type then $f:(\k^n,o)\to (\k^m,o)$ is of weakly-finite singularity type, see \S\ref{Sec.Weak.Fin.Higher.Crit.Loc.Properties}.iii.
\\
Already this gives plenty of morphisms of weakly-finite singularity type but not of finite singularity type.

   \item Take a map $f:X\to (\k^2,o),$ with $\De\sset (\k^2,o)$ a curve-germ. Then $Crit_2$ lies inside a subgerm of the germ of points at which the map
    $Crit\to \De$ is not a locally trivial fibration, see \S\ref{Sec.Weak.Fin.Critical.Locus.Properties}. If that latter germ is a point, then $f$ is of weakly-finite singularity type.

\item Suppose $Crit_X f=X,$ see Example \ref{Ex.Critical.Locus}.iv and Example \ref{Ex.Critical.Locus}.v.
 Then $X=Crit_1=Crit_2=\cdots.$ In particular, $f$ is not of weakly-finite singularity type.
     \eee
\eex

\section{\LRAP\ for analytic morphisms of weakly-finite singularity type}\label{Sec.LRAP.for.k{x}}
 Let $f:X\to Y$ be an analytic dominant morphism of weakly-finite singularity type. Take its higher critical loci,
  $X=Crit_0\supset Crit_1\supset\cdots,$ and their defining ideals, $I_{Crit_\bullet}\sset R_X.$
 Then the restriction $f_|:Crit_r\to \De_r$ is finite for some $r\ge0.$ When $char(\k)>0$   we assume   (for $N\gg1$):
\beq\label{Eq.Assumptions.for.LRAP.for.weakly.finite.sing.type}
 \bullet\  \text{All the  derivations in  $ I_{Crit_1}^N\cdot Der_{Crit_0}$,
 $I_{Crit_2}^N\cdot Der_{Crit_0,Crit_1},$\dots, $ I_{Crit_r}^N\cdot Der_{ Crit_0,\dots,Crit_{r-1}},$ }
\eeq
\hspace{1.6cm} are approximately integrable,   cf.  \S\ref{Sec.Weak.Fin.Integrability.of.Derivations}.

\hspace{0.8cm}$\bullet$
 All the automorphisms of $\quots{R_X}{I_{Crit_r}^N}$ are approximately liftable to automorphisms

 \hspace{1.2cm}of
  $\quots{R_X}{I_{Crit_{r-1}}},$ \dots, $\quots{R_X}{I_{Crit_1}},$ $\quots{R_X}{I_{Crit_0}}=R_X,$
  cf. \S\ref{Sec.Weak.Fin.Lifting.of.Automorphisms}.

\bthe\label{Thm.LRAP.weakly.finite.analytic}
 The property \LRAP (with $\Phi_X,\!\Phi_Y$ isomorphisms) holds   for   $f .$
\ethe
\bpr
  Suppose $\hPhi_Y\circ \tf=f\circ\hPhi_X,$ for some formal isomorphisms $\hPhi_X:\htX\to \hX,$  $\hPhi_Y:\htY\to \hY.$\vspace{-0.1cm}

\bee[\!\!\bf Step 1.]
\item
The formal isomorphism $\hPhi_X\!:\htX\!\to \!\hX $ sends the chain of (formal germs of) higher criti-\\cal loci,
  $\tX\!=\!Crit_0(\tf)\!\supseteq \!Crit_1(\tf) \!\supseteq\cdots,$ to the chain
$X\!=\!Crit_0(f)\!\supseteq \!Crit_1(f)\!\supseteq  \!\cdots,$ see \S\ref{Sec.Weak.Fin.Critical.Locus.Properties}.
 As before, we take the loci $Crit_\bullet,$ $\De_\bullet$ with reduced scheme structure.

Therefore \!$\hPhi_X$ \!is  a  formal solution  to \!the conditions $\Phi_X(\!\sqrt{\!Fitt_0(\cC_{j,f})})\!=\!\!\sqrt{\!Fitt_0(\cC_{j,\tf}) }\!\sset\! R_\tX,$  for $j\!\ge\! 0.$
    These are    implicit-function equations, see Remark \ref{Rem.Background.Embedding.as.IF.eq}.
 Applying the ordinary Artin approximation, Example \ref{Ex.RAP}, we   identify  $X\!=\! \tX$ and assume:
  $Crit_j(f)\!=\!Crit_j(\tf)\! \sseteq\! X$ for each $j\!\ge\!1.$
  Below we denote these subgerms by $Crit_\bullet\!\sseteq \!X.$

\medskip

Similarly, the formal isomorphism $\hPhi_Y:\htY\to \hY $ sends the chain of (formal germs of) higher discriminants,
  $\tY=\De_0(\tf)\supseteq \De_1(\tf) \supseteq\cdots,$ to the chain
$Y=\De_0(f)\supseteq \De_1(f) \supseteq\cdots.$
By the same argument as before we identify  $Y= \tY$ and assume:
  $\De_j(f)=\De_j(\tf) \sseteq Y$ for each $j\ge1.$
  Below we denote these subgerms by $\De_\bullet\sseteq Y.$

All the further automorphisms of $X,Y$ below will preserve these subgerms $Crit_\bullet,\De_\bullet$.
\item
  The restrictions $f_|,\tf_|: Crit_r\to \De_r$ are formally $\cL\cR$-equivalent and finite. Thus $f_|,\tf_|$ are  equivalent by an element of  $Aut_{Crit_r}\times Aut_{\De_r}$, by \S\ref{Sec.LRAP.finite.maps}.
 More precisely, the formal equivalence is approximated (to any order) by analytic equivalences.

These (analytic) automorphisms of $Crit_r,\De_r$ do not necessarily extend to automorphisms of the ambient spaces $X,Y.$
 Therefore, we replace the germs  $Crit_r,\De_r$ by their $N$-th infinitesimal neighborhoods, $Spec(\quots{R_X}{I_{Crit_r}^N})$ and
 $Spec(\quots{R_Y}{I_{\De_r}^N}),$ for $N\gg1.$ Here $N\gg1,$ the precise condition is given in Step 3.
  The finite maps $Spec(\quots{R_X}{I_{Crit_r}^N})\stackrel{f_|,\tf_|}{\to}Spec(\quots{R_Y}{I_{\De_r}^N})$ are still formally left-right equivalent.
  Invoke \S\ref{Sec.LRAP.finite.maps} to approximate this formal equivalence by analytic equivalence.
  The so obtained automorphisms of $Spec(\quots{R_X}{I_{Crit_r}^N})$ and  $Spec(\quots{R_Y}{I_{\De_r}^N})$ admit approximate lifting to all the groups  $Aut_{Crit_{r-1}}\times Aut_{\De_{r-1}},$ \dots, $Aut_X\times Aut_Y.$
  \bei
  \item In $char>0$ this holds by the assumption \eqref{Eq.Assumptions.for.LRAP.for.weakly.finite.sing.type}.
\item    In zero characteristic one invokes Example \ref{Ex.Lifting.Automorphisms} to get an automorphism $\Phi_X\circlearrowright R_X$
 that preserves the chain $I_{Crit_\bullet},$ and similarly $\Phi_Y\circlearrowright R_Y$ that preserves the chain  $I_{\De_\bullet}.$
\eei

Altogether, we approximate  $\hPhi_X,\hPhi_Y$ by analytic automorphisms $\Phi_X,\Phi_Y$ that preserve all the subgerms $Crit_\bullet,\De_\bullet,$ and satisfy: $\Phi_Y\circ\tf_|\circ\Phi_X^{-1}=f_|.$

 Replace $\tf$ by $\Phi_Y\circ\tf\circ\Phi^{-1}_X.$
 Thus below we assume $f_|=\tf_|:Crit_r\to \De_r,$ i.e. $f-\tf\in I_{Crit_r}^{N'}\cdot \RmX,$ where $N'\gg1$ (because $N\gg1$).
\\
For $r\!=\!0,$ with $Crit_0\!=\!X,$ i.e. $f$ is finite, this  ends the proof. Below we assume $r\!\ge\!1.$


\item
If $r=1,$ i.e. $f$ is of finite singularity type, then we invoke the finite $\cR$-determinacy by the critical locus, Lemma \ref{Thm.Critical.Locus.Fin.Determ.Map}.
 We get: $\tf=f\circ \Phi_X,$ for an analytic automorphism $\Phi_X\in Aut_X$ satisfying: $\Phi(x)-x\in I_{Crit_r}^{N''}$ for some  $N''\gg1.$

For $r\ge2$ one applies the   finite $\cR$-determinacy by higher critical loci, Lemma \ref{Thm.Higher.Critical.Locus.Fin.Determ.Map}, to get the same conclusion.
\eee

Finally, a lower bound on $N $ comes via the condition on $N',$ which is: $N'\gg1$ such that lemmas \ref{Thm.Critical.Locus.Fin.Determ.Map} and \ref{Thm.Higher.Critical.Locus.Fin.Determ.Map} apply, and the assumptions \eqref{Eq.Assumptions.for.LRAP.for.weakly.finite.sing.type} hold.
\epr

\beR
Unlike theorems \ref{Thm.LRAP.for.Nash}, \ref{Thm.LRAP.finite.analytic.maps}, here we assume that $\Phi_X,\Phi_Y$ are isomorphisms. This is used in the proof (Step 1). It is not clear, whether \LRAP holds when $\Phi_X,\Phi_Y$ are not invertible.
\eeR

  \end{document}